\documentclass{amsart}

\usepackage{amssymb}

\def\bc{\begin{center}}
\def\ec{\end{center}}
\def\be{\begin{equation}}
\def\ee{\end{equation}}
\def\ben{\begin{enumerate}}
\def\een{\end{enumerate}}
\def\bfg{\begin{figure}}
\def\efg{\end{figure}}
\def\bq{\begin{quote}}
\def\eq{\end{quote}}
\def\bd{\begin{description}}
\def\ed{\end{description}}

\def\h{\hbar}
\def\p{\partial}
\def\w{\wedge}

\def\dim{\operatorname{dim}}

\newcommand{\CC}{{\mathbb C}}
\newcommand{\RR}{{\mathbb R}}

\newcommand{\ZZ}{{\mathbb Z}}

\newcommand{\lan}{\langle}
\newcommand{\ran}{\rangle}
\newcommand{\Gb}{{\mathfrak B}}

\newcommand{\gb}{\beta}

\renewcommand{\ge}{\varepsilon}

\newcommand{\gl}{\lambda}

\newcommand{\go}{\omega}

\newcommand{\gG}{\Gamma}

\newcommand{\gD}{\Delta}

\newcommand{\gO}{\Omega}

\newcommand{\J}{{\mathcal J}}

\newcommand{\M}{\overline{\mathcal M}}

\renewcommand{\CC}{\mathbb C}
\renewcommand{\t}{\mathbf t}
\newcommand{\q}{\mathbf q}
\newcommand{\f}{\mathbf f}
\newcommand{\g}{\mathbf g}
\newcommand{\m}{\mathbf m}
\newcommand{\x}{\mathbf x}
\newcommand{\y}{\mathbf y}

\newcommand{\Q}{\mathbf Q}

\newcommand{\D}{\mathcal D}

\newcommand{\A}{\mathcal A}
\newcommand{\W}{\mathcal W}
\renewcommand{\a}{\alpha}
\renewcommand{\b}{\beta}
\renewcommand{\c}{\gamma}
\renewcommand{\d}{\delta}
\renewcommand{\H}{{\mathcal H}}
\newcommand{\F}{{\mathcal F}}
\newcommand{\G}{{\mathcal G}}
\newcommand{\T}{{\mathcal T}}
\newcommand{\C}{{\mathcal C}}
\renewcommand{\L}{{\mathcal L}}
\renewcommand{\O}{{\Omega}}
\newcommand{\1}{{\bf 1}}

\renewcommand{\l}{{\lambda}}

\newcommand{\0}{{\mathbf 0}}

\begin{document}

\title{ $A_{n-1}$ singularities and $n$KdV hierarchies }
\author{Alexander  Givental} 
\address{UC Berkeley} 
\thanks{Research is partially supported by NSF Grant DMS-0072658} 


\begin{abstract}
According to a conjecture of E. Witten \cite{W} proved by M. Kontsevich
\cite{Ko}, a certain generating function for intersection indices on the 
Deligne -- Mumford moduli spaces of Riemann surfaces coincides with a
certain tau-function of the KdV hierarchy.  
The generating function is naturally generalized under the name 
the {\em total descendent potential} in the theory of 
Gromov -- Witten invariants of symplectic manifolds.
The papers \cite{GiQ, DZ} contain two equivalent constructions,
motivated by some results in Gromov -- Witten theory, which 
associate a total descendent potential to any semisimple Frobenius
structure. In this paper, we prove that in the case of K.Saito's
Frobenius structure \cite{S} on the miniversal deformation of 
the $A_{n-1}$-singularity, the total descendent potential is  
a tau-function of the $n$KdV hierarchy. We derive this result
from a more general construction for solutions of the $n$KdV hierarchy
from $n-1$ solutions of the KdV hierarchy. 
\end{abstract}

\maketitle

\section{Introduction: Singularities and Frobenius structures.}

First examples of Frobenius structures were discovered by K. Saito \cite{S}
in the context of singularity theory. We begin with a brief
overview of very few basic elements of his (rather sophisticated) 
construction and refer to \cite{H} for further details.

Let $f: \CC^m,0\ \to\ \CC, 0 $ be the germ of a holomorphic function 
at an isolated critical point of multiplicity $N$. 
We will assume for simplicity that $f$ is
weighted-homogeneous. Let $\T $ be the parameter space of its
miniversal deformation $F(x,\tau)$. Tangent spaces to $\T$ are naturally
equipped with the algebra structure: 
$T_{\tau} \T = \CC [x] /(F_x(\cdot,\tau))$.
Pick a holomorphic weighted-homogeneous volume form $\go_{\tau}$ on $\CC^m$
possibly depending on the parameters $\tau$. Then the Hessians $\Delta (x)$
of critical points $x\in \operatorname{crit}( F(\cdot, \tau)) $
become well-defined. The corresponding residue paring 
\[ (\phi , \psi )_{\tau} = \sum_{x\in \operatorname{crit}( F(\cdot, \tau))} 
\frac{\phi(x)\psi(x)}{\Delta (x)} \]
is known to define a non-degenerate symmetric bilinear form on 
$T_{\tau}\T$ which depends analytically on $\tau$, 
extends across the bifurcation hypersurface without singularities
and thus makes $T_{\tau}\T$ Frobenius algebras. The key point of K.Saito's
theory is that there exists (according to a theorem of M. Saito, see \cite{H})
a choice of $\go$ (called {\em primitive}) 
that makes the family of Frobenius algebras a Frobenius structure. 
The latter means certain integrability property
which will be recalled lated when needed.
We refer to \cite{D, Ma} for a detailed account of numerous manifestations 
of the property --- such as flatness of the metric $(\cdot,\cdot)$ for
example.

In the case of simple singularities a weighted-homogeneous volume
form $\go $ coincides with $dx_1\w ... \w dx_m$ 
(up to a non-zero constant factor
which in fact does not affect the metric $(\cdot,\cdot)$ ) 
and therefore $\go $ is primitive. 

In the example $A_{n-1}$ we set
$m=1$, $f(x)=x^n/n $, $F(x,\tau)=x^n/n+\tau_1x^{n-2}+...+\tau_{n-1}$,
$\go = dx$. The basis $\{ \p_{\tau_i} \} $ in $T_0\T$ is identified
with the basis $x^{n-1-i}$ of the local algebra $H=\CC [x]/(x^{n-1})$,
and the residue pairing in this basis takes the form 
$(\p_{\tau_i},\p_{\tau_j})_0=\delta_{i+j,n-1}$. By the general theory,
the following residue metric is flat:
\[ (\p_{\tau_i},\p_{\tau_j})_{\tau} = 
\operatorname{Res}_{x=\infty}
\frac{x^{2n-2-i-j}dx}{F'(x,\tau)}.\]

\medskip

In Gromov -- Witten theory, intersection indices in moduli spaces of genus-$0$ 
pseudo-holomorphic curves in a given compact symplectic manifold define a Frobenius
structure on the cohomology space of the manifold. What is the structure behind
intersection theory in spaces of higher genus pseudo-holomorphic curves, and is 
it possible to recover the totality of higher genus Gromov -- Witten invariants
from the Frobenius structure? While the answer to the first question is yet
unknown, the answer to the second one seems to be positive in the semisimple 
case.

According to \cite{GiQ} {\em the total descendent potential} corresponding to
a semisimple Frobenius manifold can be defined by the formula
\begin{equation} \label{1.1}
\D (\q)= C(\tau)\ \hat{S}_{\tau}^{-1}\ \Psi (\tau)\ \hat{R}_{\tau}\ 
\exp (U/z)\hat{\ }\ \prod_{i=1}^N \D_{A_1}(\q_i). \end{equation}
The ingredients of the formula will be explained 
later in the context of singularity theory. 
Roughly, the function $\ln \D$ is supposed to have the form 
of ``a genus expansion''
$\sum \h^{g-1} \F^{(g)}(\q)$ where $\F^{(g)}$ depend
on the sequence $\q$ of vector variables $q_0,q_1,q_2,... $ taking values
in the local algebra $H$ of the singularity.
The Taylor coefficients of $\F^{(g)}$ are to play the role of 
genus-$g$ Gromov -- Witten invariants and their gravitational descendents.
The product term in 
(\ref{1.1}) is the tensor product of
N copies of the total descendent potential for the $A_1$-singularity
(which is a tau-function of the KdV hierarchy and is discussed in
Section 3). The product is considered as an ``element of a Fock space''.
The $S$, $R$ and $\exp (U/z)$ are elements of a certain group (of loops
in the variable $z$) acting on the 
elements of the Fock space via some ``quantization'' representation 
$\hat{\ }$. The loops $S(z)$, $R(z)$ and $\exp (U/z)$ 
(as well as $C$ and $\Psi $ which 
are a non-zero normalizing constant and an invertible matrix)
are defined in terms of the Frobenius structure and in the case 
of singularities allow convenient descriptions via oscillatory integrals and
their asymptotics. The ingredients of the formula depend on
a choice of the point $\tau \in \T$ which has to be semisimple,  
i. e. the function $F(\cdot , \tau)$ must have $N$ non-degenerate critical
points. For example, $U$ is the diagonal matrix of the critical values
of $F(\cdot, \tau)$. As it is explained in \cite{GiQ},  
the resulting function $\D $ does not depend on $\tau$, satisfies the
so called {\em $3g-2$-jet condition}, {\em Virasoro constraints} and 
has the correct (in the sense of \cite{D}) genus-$0$ part $\F^{(0)}$.
\footnote{According to a result from \cite{DZ}, a function with these 
properties, when exists, is unique.} 

In this paper, we will prove that in the case of $A_{n-1}$-singularities,
the function (\ref{1.1}) is a tau-function of the $n$KdV-hierarchy (Theorem $1$).
 
In Section $2$, we describe the quantization formalism underlying (\ref{1.1}).
The KP, KdV and nKdV-hierarchies are described in Section $3$ in terms of
the so called {\em vertex operators} of the infinite dimensional Lie
algebra theory \cite{K}. In Section $4$, we reconcile the notations 
of representation theory and singularity theory and state Theorem $1$.
In Sections $5$ and $6$, we study conjugations of the vertex operators
by the operators $S$ and $R$. The corresponding Theorems $2$ and $3$ 
are the technical heart of the paper and provide surprisingly simple and 
general formulations in terms of singularity theory. In Section $7$, 
we show how various central constants (somewhat neglected in the 
preceeding computations) are governed by a certain multiple-valued closed 
$1$-form $\W$ on the complement to the discriminant.
The form $\W$ appears to be a new object in singularity theory, and its 
properties play a key role in the proof of Theorem $1$. 
In Section $8$, we discuss in detail the ``Fock spaces'' intertwined by 
the operators $S$ and $R$ and describe analyticity properties of
our vertex operators. In Section $9$, we state and prove Theorem $4$ 
which interprets the formula (\ref{1.1}) as a device transforming some 
solutions of the KdV-hierarchy into solutions of the $n$KdV-hierarchy 
(and which contains Theorem $1$ as a special case). Relationships with
``$W_n$-gravity theory'' are discussed in Section $10$ (Theorem $5$).        
The appendix, included mostly for aesthetic considerations, contains a 
direct treatment of genus-$0$ consequences of Theorem $1$.

Slightly generalizing the methods of the present paper, one can prove
that the total descendent potential (\ref{1.1}) corresponding to an
ADE-singularity satisfies an integrable hierarchy described explicitly
in terms of vertex operators and very similar to the famous 
hierarchy of Kac -- Wakimoto \cite{KW} constructed via representation 
theory of loop Lie algebras. We will return to this subject in \cite{GiM}   

\medskip

{\em Acknowledgments.} Substantial part of the paper was written
during our stay at IHES (Paris) and MPI (Bonn) in Summer '$02$. 
We would like to thank these institutions for hospitality, and
the National Science Foundation --- for financial support. 
We are also thankful to
E. Frenkel and P. Pribik for their interest and stimulating discussions,
to A. Schwarz for consultations on $W_n$-gravity, 
and especially to T. Milanov for several useful observations.  

\section{The quantization formalism.}
 
Consider the local algebra $H=\CC[x]/(f_x)$ as a vector space
with a non-degenerate symmetric bilinear form defined by the  
residue pairing 
\[ ( a, b )_0 = \operatorname{Res}_{x=0}  
a(x) b(x) dx_1\w ...\w dx_m/f_{x_1}...f_{x_m}.\] 
Let $\H = H((z^{-1}))$ denote the space 
of Laurent series in one indeterminate $z^{-1}$ with coefficients in $H$.
We equip $\H$ with the even symplectic form
\begin{equation} \label{2.1} 
 \O (\f,\g) =\frac{1}{2\pi i} \oint ( \f(-z), \g(z))_0 dz =-\O (\g, \f). 
\end{equation}    
The polarization $\H=\H_{+}\oplus \H_{-}$ defined by the 
lagrangian subspaces $\H_{+}=H [z]$, $\H_{-}=z^{-1} H [[z^{-1}]]$ identifies
$(\H, \O)$ with the cotangent bundle $T^*\H_{+}$. Then the standard
quantization convention associates to constant, linear and 
quadratic hamiltonians $G$ on $(\H, \O)$ 
differential operators $\hat{G}$ of order $\leq 2$ acting on functions on 
$\H_{+}$. More precisely, let $\{ q_{\a} \}$ be a coordinate 
system on $\H_{+}$ and $\{ p_{\a} \}$ --- the dual coordinate system on 
$\H_{-}$ so that the symplectic structure in these coordinates assumes
the Darboux form $\O  = \sum_{\a} p_{\a} \w q_{\a}$. 
For example, when $H$ is the standard one-dimensional Euclidean space then
\begin{equation} \label{2.4}
\f= \sum q_k z^k + \sum p_k (-z)^{-1-k} \end{equation}
is such a coordinate system.
In a Darboux coordinate system the quantization convention reads
\begin{equation} \label{2.2}
 (q_{\a})\hat{\ }:=  q_{\a}/\sqrt{\h}, \ \ (p_{\a})\hat{\ } := 
\sqrt{\h}\ \p /\p q_{\a} ,\end{equation}   
\begin{equation} \label{2.3}
 (q_{\a}q_{\b})\hat{\ } := \frac{q_{\a}q_{\b}}{\h} ,\ \ 
(q_{\a}p_{\b})\hat{\ } := q_{\a} \frac{\p}{\p q_{\b}}, \ \ 
(p_{\a}p_{\b})\hat{\ } := \h \frac{\p^2}{\p q_{\a}\p q_{\b}} .\end{equation}
The quantization is a representation of the Heisenberg algebra of constant and
linear hamiltonians, but it is only a {\em projective} representation of the 
Lie algebra of quadratic hamiltonians on $\H$ to the Lie algebra of 
differential operators.
For quadratic hamiltonians $F$ and $G$ we have
  \[ \{ F, G\}\hat{\ } = [\hat{F}, \hat{G}] + \C (F,G) \]
where $\{ \cdot , \cdot \}$ is the Poisson bracket, $[\cdot , \cdot ]$ is the
commutator, and $\C $ is a cocycle characterized by
the properties that
\[ \C (p_{\a}p_{\b}, q_{\a}q_{\b}) = 1 \ 
\text{if}\ \a \neq \b,\ \ \C (p_{\a}^2, q_{\a}^2) = 2 \ , \]
and $\C =0$ on all other pairs of quadratic Darboux monomials.  

The differential operators act on formal functions 
(with coefficients depending on $\h^{\pm 1/2}$) 
on the space $\H_{+}$ of vector-polynomials 
$\q = q_0+q_1z+q_2z^2+...$  with the coefficients $q_0,q_1,q_2 ... \in H$.
We will often refer to such functions as {\em elements of the Fock space}. 

Consider now linear operators on $\H $ which preserve the symplectic 
structure and commute with multiplication by $z$. They form a twisted
version of the loop group $\L GL(H)$. It consists of the loops
$M(z)$ satisfying $M^t(-z)M(z)=\1$ where $\ ^t$ means transposition
with respect to the inner product $(\cdot , \cdot )_0$. Quantized operators
$\hat{M}$ are defined as $\exp (\ln M)\hat{\ }$ (though the domain
of $\hat{M}$ in the ``Fock space'' may depend on $M$). The operators
$\hat{S}$ and $\hat{R}$ in the formula (\ref{1.1}) are of this nature.
Moreover, the loops $S(z)$ and $R(z)$ are triangular in the sense
that $S(z)=\1+S_1z^{-1}+S_2z^{-2}+... $ and $R(z)=\1+R_1z+R_2z^2+...$.

\section{Example: KP and KdV hierarchies.}

The goal of this section is to reconcile the conventional theory of
integrable hierarchies with the quantization formalism of the previous 
section in the example of KdV (i.e. $2$KdV) hierarchy. 
The $n$KdV hierarchies will be
treated in this paper as ``reductions modulo $n$'' of the KP hierarchy.
The KP hierarchy has an abstract description as a sequence of commuting
flows on the semi-infinite grassmannian with the time variables 
$x_1,x_2,x_3,...$. The ``bosonic-fermionic correspondence'' identifies
the space of semi-infinite forms with the symmetric algebra 
$\CC [\x]$ in the variables $\x= (x_1,x_2,x_3,... )$. 
Under the Pl\"ucker embedding, points of the grassmannian
are transformed into 1-dimensional subspaces spanned by   
certain functions of $\x$, and the KP flows are
defined tautologically as time translations. The equations of the 
KP hierarchy thus assume the form of {\em Hirota quadratic equations}
describing the image of the grassmannian under the Pl\"ucker embedding.

It will be convenient for us to use the following
{\em vertex operator construction} of the Hirota quadratic equations. 
According to \cite{K}, Ch. $14$, a function 
$\Phi (\x)$ (which we will assume to have the form
$\exp \sum \h^{g-1} \phi^{(g)}(\x) $) satisfies the KP hierarchy iff
\begin{equation} \label{3.1}
\operatorname{Res}_{\zeta =\infty} d\zeta \ 
e^{\sum_{j>0}\zeta^j(x'_j-x''_j)/\sqrt{\h}} \ 
e^{-\sum_{j>0}\frac{\zeta^{-j}}{j} \sqrt{\h}(\p_{x_j'}-\p_{x_j''})}
\ \Phi(\x')\Phi(\x'') =0. \end{equation}
The equation is interpreted in the following way. The change 
\[
x_j=(x_j'+x_j'')/2,\ y_j=(x_j'-x_j'')/2,\ \ \p_{x_j}=\p_{x_j'}+\p_{x_j''},
\ \p_{y_j}=\p_{x_j'}-\p_{x_j''} \]
transforms the equation (\ref{3.1}) into
\begin{equation} \label{3.2} 
\operatorname{Res}_{\zeta =\infty} d\zeta \ 
e^{2\sum_{j>0}\zeta^jy_j/\sqrt{\h}} \ 
e^{-\sum_{j>0}\frac{\zeta^{-j}}{j} \sqrt{\h} \p_{y_j}}
\ \Phi(\x+\y)\Phi(\x-\y) =0. \end{equation}
Expanding in $\y$ yields an infinite system of equations
on partial derivatives of $\Phi (\x)$ which is an abstract
form of the KP hierarchy.

Note that prior to extracting the residue, the expansion of (\ref{3.2}) in 
$\y$ is an infinite series with the property
that the coefficient at each monomial $\y^{\bf m}$ is a Laurent series
in $\zeta^{-1}$, i.e. the powers of $\zeta$ are bounded from above
by a constant depending on $m$. We should therefore think of the
expressions in (\ref{3.1}),(\ref{3.2}) as expansions near $\zeta =\infty$.
Below we call such an expression {\em regular in $\zeta$} if it 
contains no negative powers of $\zeta $, i.e. 
the coefficient at each monomial $\y^{\bf m}$ is a polynomial.

By definition, solutions of the {\em $n$KdV hierarchy} (also called  
{\em Gelfand -- Dickey} or {\em $W_n$-hierarchy}) are those solutions
of the KP hierarchy which do not depend on $x_j$ with 
$j \equiv 0\ \operatorname{mod}\ n$. For $n=2$ we obtain the
KdV hierarchy whose solutions depend therefore only on $\x_{\text{odd}}$ and
do not depend on $\x_{\text{even}}$. Note that the derivations
$\p_{y_{2k}}$ in (\ref{3.2}) can be omitted while the multiplications by
$y_{2k}$ cannot. Thinking of $\exp 2\sum_{k>0} \zeta^{2k} y_{2k}$ as 
an arbitrary function of $\zeta^2$ and symmetrizing (\ref{3.2})
over the Galois group $\ZZ_2$ of the covering $\zeta \mapsto \zeta^2$,
we arrive at the following description of the KdV hierarchy:

{\em a function $\Phi (\x_{\text{odd}})$ satisfies the KdV hierarchy 
if and only if the following differential $1$-form  
is regular in $\zeta^2$:}
\begin{equation} \label{3.3} 
\sum_{\pm} \pm d \zeta \ 
e^{\pm \sum_{j\ \text{odd}}\zeta^j(x'_j-x''_j)/\sqrt{\h}} \ 
e^{\mp \sum_{j\ \text{odd}}\frac{\zeta^{-j}}{j} \sqrt{\h}
(\p_{x_j'}-\p_{x_j''})} \ \Phi(\x')\Phi(\x'') .\end{equation}

\medskip

The {\em Witten-Kontsevich tau-function} is defined as
\begin{equation} \label{3.4} 
\T (\t)=\exp \sum_{g=0}^{\infty} \h^{g-1} \sum_{m=0}^{\infty}\frac{1}{m!}
\int_{\M_{g,m}} \t(\psi_1)\w ... \w \t(\psi_m) ,\end{equation} 
where $\M_{g,m}$ are the Deligne -- Mumford moduli spaces of stable genus
$g$ compact complex curves with $m$ marked points, $\psi_i$ are 
the $1$-st Chern classes of the {\em universal cotangent line} bundles
(formed by the cotangent lines to the curves at the $i$-th marked points)
over $\M_{g,m}$, and $\t$ is a polynomial $\t (z) = t_0+t_1 z+ t_2 z^2 + ...$.
It is known (see for instance \cite{W}) that $\T$ satisfies the 
{\em string} and {\em dilaton} equations
\[ \p_{t_0} \T - \sum_{k=0}^{\infty} t_{k+1}\p_{t_k}\T =\frac{t_0^2}{2\h}\T,\  
\ 3\p_{t_1}\T-\sum_{k=0}^{\infty} (2k+1) t_k\p_{t_k}\T = -\frac{1}{8}\T .\]
At $t_0=0$ the genus-$g$ part of $\ln \T$ depends only on $t_1,...,t_{3g-2}$
(for dimensional reasons). This implies that $\T$ is well-defined at least
as a formal function of, say, $\h, t_0/\h, t_1, t_2, ...$. 
Note however that the vector fields on the  
LHS of the string and dilaton equations become linear homogeneous
after the change of variables $q_k=t_k-\d_{k,1}$ 
called the {\em dilaton shift}. We define an element in the Fock space by
\begin{equation} \label{3.5} D_{A_1}(\q ) := \T ( \t ),\ \text{where}
\ \q(z):=\t(z)-z .\end{equation}
Thus $\D_{A_1}$ is well-defined as a formal function near the shifted 
origin $q(z)=-z$. 

\medskip

According to Witten's conjecture \cite{W} proved by Kontsevich \cite{Ko}
the function $\D_{A_1}$ satisfies the KdV hierarchy (\ref{3.3}) after the
substitution $q_k = (2k+1)!! x_{2k+1}, \ k=0,1,2,...$. We also have
$\p_{x_{2k+1}}= (2k+1)!! \p_{q_k}$. We are going to rewrite
(\ref{3.3}) in terms of Section 2.  The exponents in (\ref{3.3}) are
elements of the Heisenberg Lie algebra and are quantizations of linear 
hamiltonians in the symplectic space $\H$. We will encode the hamiltonians
by the corresponding (constant) hamiltonian vector fields. 
The standard relationship
$\dot{q}=h_p, \ \dot{p}=-h_q$ between hamiltonians $h$ and their vector fields
dictates the following correspondence 
between the Darboux coordinates (\ref{2.4})
as linear functions on $\H$ and vectors in $\H$:
$p_k \mapsto  z^k,\ \ q_k \mapsto -(-z)^{-1-k}, \ \ k=0,1,2,... $.
Using the notation $\gl=\zeta^2/2$, we can rewrite the KdV hierarchy 
(\ref{3.3}) for $\D_{A_1}$ in the form    
\begin{equation} \label{3.6}   \sum_{\sqrt{2\gl}=\pm \zeta} 
\ (\gG^{-}(\gl) \D_{A_1})(\q')\ \  (\gG^{+} (\gl) \D_{A_1})(\q'')\ \ 
\frac{d \gl}{\sqrt{\gl}}
 \ \ \text{is regular in $\gl$},\end{equation}
where the sum is taken over the two values of $\sqrt{2\gl}$, and 
\begin{equation} \label{3.7}
\gG^{\pm}(\gl):= 
e^{\pm \sum_{k<0} (\frac{d}{d\gl})^k (2\gl)^{-1/2}\ (-z)^k }\ 
e^{\pm \sum_{k\geq 0} (\frac{d}{d\gl})^k (2\gl)^{-1/2}\ (-z)^k } . 
\end{equation}
We will informally refer to (\ref{3.6},\ref{3.7}) as 
the KdV hierarchy for the total descendent potential $\D_{A_1}$.

\section{The vertex operators for $n$KdV.}

Returning to the setting of Section $1$, we introduce vertex operators 
associated with cycles vanishing at isolated critical points in a fashion
generalizing the role of
\[ (2\gl)^{-1/2}=\int_{[x]: x^2/2=\gl} dx/d(x^2/2)  \]
in (\ref{3.7}). More precisely, the operators will have the form
\begin{equation} \label{4.1}
\gG^{\beta} =e^{  \sum_{k<0} I_{\beta}^{(k)}(\gl)(-z)^k }\  
e^{\sum_{k\geq 0} I_{\beta}^{(k)}(\gl) (-z)^k} , \end{equation}
where $I_{\beta}^{(k)}$ are vector functions with values in $H$ 
which are consecutive derivatives of one another,
$d I^{(k)}_{\beta} /d\gl = I^{(k+1)}_{\beta}$, and are defined as
follows.    
 
Let $f$ be a weighted-homogeneous singularity with the local algebra $H$
and with the residue pairing $(\cdot,\cdot)_0$ defined by the volume form 
$\go_0=dx_1\w ... \w dx_m$. We will always assume that the number of variables
$m=2l+1$ is odd, that the monomials $\phi_1,...,\phi_N=1 \in \CC [x]$ 
represent a basis in $H$, 
and that the spectrum $\deg (\phi_1 \go),...,\deg (\phi_N\go)$ 
contains no integers.
For $[\phi] \in H$ represented by a linear combination $\phi $ of the 
monomials $\phi_i$, we put 
\[ (\ I^{(0)}_{\beta}(\gl), \ [\phi] \ )_0 := 
(\frac{1}{2\pi}\frac{d}{d\gl})^l\ \int_{\beta \subset f^{-1}(\gl)}
\phi (x) \frac{dx_1\w ... \w dx_m}{df (x)} ,\]
where $\beta $ is a middle-dimensional cycle in the Milnor fiber $f^{-1}(\gl)$.
\footnote{ When $\gb$ is the vanishing cycle $f^{-1}(\gl)\cap \RR^m$ of the
$A_1$ singularity $f=(x_1^2+...+x_m^2)/2$ in $m=2l+1$ variables, 
we have $\int_{\beta} dx / df =
\sigma_{2l} \gl^{l-1/2}$ where $\sigma_{2l} =2 (2\pi)^l/(2l-1)!!$ is the 
volume of the unit $2l$-dimensional sphere. The factor $1/(2\pi)^{l}$ in the
definition of $I^{(0)}$ 
makes therefore $I_{\beta}^{(0)}=2\gl^{-1/2}$ independent on $l$.} 
This defines $I_{\beta}^{(0)}(\gl)$ as a vector-function with homogeneous
components of non-integer degrees, and we extend the definition to 
$I^{(k)}_{\beta}$ by the obvious derivations and anti-derivations in $\gl$. 
This determines the vertex operator unambiguously up to the classical 
monodromy of the cycle. 

\medskip

In the special case of $A_{n-1}$ singularities, $\gG^{\beta}$
are closely related to the vertex operators of the $n$KdV hierarchy.
Put $l=0$, $f=x^n/n$, $\phi_i=x^{n-1-i}, \ i=1,...,n-1$. 
We take the cycle $\beta$ to
be {\em one point} $x=(n\gl)^{1/n}$ at the level $f^{-1}(\gl)$ and denote
this cycle $\a$. Then  
\[ (I_{\a}^{(0)}, [\phi_i ])_0= 
\int_{\a} x^{n-1-i} \frac{dx}{d x^n/n}  = (n\gl)^{-i/n} .\]  
Equivalently, $I_{\a}^{(0)} = \sum_{i=1}^{n-1} [x^{i-1}] (n\gl)^{-i/n}$.
This implies
\[ \sum_{k\in \ZZ} I_{\a}^{(k)}(-z)^k=\sum_{i=1}^{n-1} \sum_{k\in \ZZ} 
[x^{i-1}] z^k\ (n\gl)^{-(i+kn)/n} 
\prod_{r=0}^{\infty} (i+rn) / \prod_{r=k}^{\infty} (i+rn) .\]
The double sum contains exactly one summand with each power $i+kn$ of 
$\zeta = (n\gl)^{1/n}$ not divisible by $n$. 

Let us compare the coefficients at $\zeta^{-j}$ and $\zeta^j$.
For $j=i+kn$ we have $-j=n-i+(-1-k)n$. The corresponding vectors 
$[x^{i-1}]z^k$ and $[x^{n-i-1}]z^{-1-k}$ in $\H = H((z^{-1}))$
have the symplectic inner product $(-1)^k$ (while any other pairs are
$\gO$-orthogonal). The corresponding 
factorial products multiply to $(-1)^{k+1}/(i+kn)$. 

Let $\p/\p q_{i,k}$ denote the elements in the Heisenberg algebra (acting
on the Fock space of functions on $\H_{+}$) which correspond to 
the vectors $[x^{i-1}] z^k$ in $\H$. The above computation means that
the change 
\begin{equation} \label{4.2} q_{i,k} = i (i+n) (i+2n) ... (i+kn) x_{i+kn} 
\end{equation}
transforms $\sum_{k\in\ZZ} I^{(k)}_{\a} (-z)^k$ into
$ -\sum_{j < 0} x_j\ \zeta^j + 
\sum_{j > 0}\p_{x_j}\ \zeta^{-j}/j $ where $j\in\ZZ\backslash n\ZZ$. 

Comparing with (\ref{3.1}) we see that the change (\ref{4.2}) transforms
solutions of the $n$KdV hierarchy into functions $\D $ satisfying the
condition 
\begin{equation} \label{4.3} 
\sum_{\a}\ (\gG^{-\a} \D) (\q')\ (\gG^{\a} \D) (\q'') \ 
\gl^{(1-n)/n} d\gl \ \ \text{is regular in $\gl$.} 
\end {equation}    
The sum here is taken over all the $n$ values of $\l^{1/n}$ which
correspond to the one-point cycles $\a$. In particular, the coefficients $\gl^{(1-n)/n}$
in different summands differ by appropriate $n$-th roots of unity (rather than coincide).

Our goal in this paper is to prove the following theorem.

\medskip

{\bf Theorem 1.} {\em The total descendent potential $\D_{A_{n-1}}$ of the
$A_{n-1}$-singularity defined by the formula (\ref{1.1}) (as explained in 
\cite{GiQ}) satisfies (\ref{4.3}) and therefore is transformed by the change
(\ref{4.2}) into a tau-function of the $n$KdV hierarchy.}

\section{From descendents to ancestors.}

According to the definition (\ref{1.1}) the function $\D$ in Theorem $1$
has the form 
$\D = e^{{\mathbf F}^{(1)}(\tau)}\hat{S}_{\tau}^{-1} \A_{\tau}$ 
where $\A_{\tau}$ is
some other element of the Fock space depending on $\tau\in \T$ and
called in \cite{GiQ} the {\em total ancestor potential} 
(and ${\mathbf F}^{(1)}$ is a function of $\tau$ called {\em the genus-$1$
Gromov-Witten potential} which will be described in the next section
and which actually vanishes in the case of simple singularities). 
 Replacing in (\ref{4.3}) the function $\D$ with $\A=\hat{S} \D$ and 
$\gG^{\pm\a}$ ---
with its conjugation $\hat{S} \gG^{\pm\a} \hat{S}^{-1}$ we obtain a
reformulation of Theorem $1$ in terms of the ancestor potential.
Let us compute $\hat{S} \gG^{\beta} \hat{S}^{-1}$, first formally,
and then in the actual setting of singularity theory.

A quantized ``lower-triangular'' symplectic operator $S(z)=\1+S_1z^{-1}+S_2z^{-2}+...$
acts on elements of the Fock space by the formula (Proposition $5.3$
in \cite{GiQ})
\[ (\hat{S}^{-1} \G) (\q ) = e^{W(\q,\q)/2\h} \G ([S \q]_{+}), \]
where $[S\q]_{+}$ is the truncation of negative powers of $z$ in 
$S(z) \q(z)$, and the quadratic form $W(\q,\q)=\sum (W_{kl}q_k,q_l)$ 
is defined by 
\begin{equation} \label{5.1}
 \sum_{k,l\geq 0} \frac{W_{kl}}{w^kz^l}:=\frac{S^t(w)S(z)-\1}{w^{-1}+z^{-1}}.
\end{equation} 
Respectively, 
\[ (\hat{S} \G)(\q) = 
e^{-W([S^{-1}\q]_{+},[S^{-1}\q]_{+})/2\h} \G ([S^{-1}\q]_{+}) .\]
For $\f \in H[[z,z^{-1}]]$, let
$(e^{\f})\hat{\ }\ :=e^{\hat{\f}_{-}} e^{\hat{\f}_{+}}$ be the corresponding
element in the Heisenberg group.
The previous formulas show that 
\begin{equation} \label{5.2} \hat{S}\ (e^{\f})\hat{ }\ \hat{S}^{-1} \G = 
e^{W(\f_{+},\f_{+})/2} (e^{S\f})\hat{\ }\ \G .\end{equation}   

\medskip

We are returning to the Frobenius structure on the parameter space $\T$ of 
a miniversal deformation of a (weighted - homogeneous) singularity.
Consider the complex oscillating integral
\[  \J_{\Gb}(\tau) = (- 2\pi  z)^{-m/2}\int_{\Gb} e^{F(x,\tau)/z} \go .\]
Here $\Gb $ is a non-compact cycle from the relative homology group
\footnote{The present description of the oscillating 
integral is accurate only for subdeformation $\tau\in \T^{lower}$
of $f$ by terms of degrees lower than $\deg f =1$. Our excuses are that 
(i) such $\tau$ will suffice for all our goals and (ii) $\T^{lower}=\T$
for $A_{n-1}$ and other simple singularities.}  
\[ \lim_{M\to \infty}\ 
H_m (\CC^m, \{ x: \operatorname{Re} (F(x,\tau)/z) \leq -M \})\simeq \ZZ^N. \]
We will assume that $\go $ is primitive and use the notation
$\p_1,...,\p_N$ for partial derivative with respect to a flat (and 
weighted - homogeneous)
coordinate system $(t_1,...,t_N)$ of the residue metric. Saito's theory 
of primitive forms guarantees that the differential equations 
for $\J_{\Gb}$ in flat coordinates assume the following form: 
\[ z \p_i \p_j \J = \sum_k a_{ij}^k \p_k \J, \ \ \text{where} \ \ 
\p_i \bullet \p_j =\sum_k a_{ij}^k(\tau) \p_k   \]
is the multiplication on the tangent spaces $T_{\tau}\T$.
In particular, the linear pencil of connections on the {\em cotangent} bundle
\[ \nabla := d - z^{-1}\sum (\p_i\bullet )^t dt_i \]
is flat for any $z\neq 0$ (since $ \sum (\p_j \J_{\Gb}) dt_j$ provide
a basis of $\nabla$-flat sections). The integrability of $\nabla$ is
a key axiom in the definition of Frobenius structures \cite{D}.
\footnote{Note that the operators $\p_i \bullet$ are self-adjoint with respect 
to the metric. Identifying the tangent and cotangent spaces via the metric, 
we get $\nabla=d-z^{-1}\sum (\p_i \bullet) dt_i$, while the natural
adjoint connection on the tangent bundle reads 
$d + z^{-1} \sum (\p_i \bullet)\ dt_i$.}

The oscillating integral also satisfies the following homogeneity condition:
\[ (z\p_z + \sum (\deg{t_i}) t_i \p_i )\ z\p_j \J = -\mu_j\ z\p_j \J ,\]
where $-\mu_j =  \deg (\p_j F) +\deg (\go)- m/2,\ j=1,...,N$, is the 
{\em spectrum} of the singularity symmetric about $0$. One can extend 
therefore the connection $\nabla$ to the $z$-direction by 
\[ \nabla_{\p_z} := \p_z +\mu/z +(E\bullet)^t/z^2, \]
where $E=\sum (\deg{t_i}) t_i \p_i$ is the {\em Euler field} and
$\mu =\operatorname{diag}(\mu_1,...,\mu_N)$ is the {\em Hodge 
grading operator} (anti-symmetric with respect to the metric and 
diagonal in a graded basis). The extended connection is flat (since
$\sum z\p_i \J_{\Gb}\ dt_i$ provide a basis of flat sections)
and can be considered as an {\em isomonodromic} family of connections
in $z\in \CC\backslash 0$ depending on the parameter $\tau\in \T$. 
Identifying the $T^*\T$ with $T\T$ via the metric, we obtain
the connection operator $\p_z-\mu/z+(E\bullet )/z^2$.
The connection is regular at $z=\infty$.
At $\tau =0$ it turns into $\p_z - \mu /z$.           

\medskip

{\em Definition.}
The operator $S_{\tau}(z)=\1+S_1z^{-1}+S_2z^{-2}+...$ is defined
as a gauge transformation in the twisted loop group (i.e. $S^t(-z)S(z)=\1$)
which transforms near $z=\infty $ the
connection operator $\p_z-\mu/z+(E\bullet )/z^2$ at the parameter value
$\tau\in \T$ into the connection operator $\p_z -\mu /z$.

\medskip

In particular, the basis of flat sections for the extended connection
defined by the complex oscillating integrals $z\p_i\J_{\Gb}$  near $z=\infty$
has the form $S_{\tau}(z) z^{\mu} C$ where $C$ is a constant invertible 
matrix depending on the basis of cycles $\Gb$. 
This implies that $S$ is a fundamental solution to 
$z\p_i S = \p_i\bullet S,\ i=1,...,N$, and
satisfies the homogeneity condition $(z\p_z +E) S = \mu S - S \mu$.
A choice of the series solution $S$ with these properties and satisfying 
the asymptotical condition $S(\infty) =\1$ and the symplectic condition
$S^t(-z)S(z)=\1$ is called in \cite{GiQ}
{\em calibration} of the corresponding Frobenius structure.
In general calibration is not unique (and may depend on finitely many 
constants) unless there is no integers among
the spectral differences $\mu_i-\mu_j$. 
It is therefore unique in the case of simple singularities.
For more general weighted - homogeneous singularities a canonical
choice is specified by the condition that $S_{\tau=0}=\1$. 
\footnote{This makes the descendent potential a special case of the ancestor potential
$\A_{\tau}=\hat{S}_{\tau} \D$ with $\tau=0$ provided that $\hat{S}_{0} \D$ is well-defined
(see Section $9$).} 

\medskip

Let us consider now {\em period vectors} $I_{\beta}^{(k)}$ of $\gl$ and $\tau$
with values in $H$ defined by the integrals over vanishing cycles 
$\beta \in H_{m-1}(V_{\tau}(\gl ))$ in the Milnor fibers
$V_{\gl,\tau} = \{x\in \CC^n :  F (x, \tau)=\gl \}$. We keep the
notation $m=2l+1$ and other hypotheses of Section $4$ and define the period vectors
$I^{(k)}_{\beta}$ by
\footnote{Here $d^{-1}\go$ is any $m-1$-form whose differential in $x$ equals $\go$. In the
second equality we assume for simplicity that $\go $ is independent of $\tau$ as 
is the case for simple singularities.}
 \[ ( I^{(k)}_{\beta} (\gl, \tau), \p_i ) := - (2\pi)^{-l} 
\p_{\gl}^{l+k} \p_i\int_{\beta \subset V_{\gl, \tau}} d^{-1}\go = 
(2\pi)^{-l} \p_{\gl}^{l+k} \int_{\beta} (\p_i F)\frac{\go}{dF}.\]
The vector-valued functions $I^{(k)}_{\beta}$ are multiple-valued and
ramified along the discriminant where $V_{\gl,\tau}$ becomes singular.
We refer to \cite{AVG} for a standard description of the reflection 
monodromy group for the cycle $\beta$ and the integrals.
When $\tau =0$, the vector-functions $I^{(k)}_{\beta}$ specialize
to those of the previous section.

\medskip

\noindent {\bf Theorem 2.} {\em Let $\f (\gl,\tau) = 
\sum_{k\in \ZZ} I^{(k)}_{\beta}(\gl, \tau) (-z)^k $.
Then $ \f (\gl, \tau) = S_{\tau}(z) \f (\gl, 0)$.}

\medskip

{\em Remark.} The integrals $I^{(k)}(\gl)$ expand near $\gl=\infty$ into
Laurent series (with fractional exponents), and the maximal exponent in 
$I^{(k)}$ tends to $-\infty $ as $k\to \infty$. Respectively, 
coefficients in a $z$-series 
of the form $S \f$ with $S=\sum_{l\geq 0} S_l z^{-l}$ and 
$\f=\sum I^{(k)}(-z)^k$, which are infinite sums 
$\sum_{l\geq 0} \pm S_l I^{(k+l)}$, converge in the $1/\gl$-adic sense. 

\medskip
  
{\em Proof.}
The period vectors $I^{(k)}_{\beta}$ are related to the oscillating 
integrals $J_{\Gb}$ by (a version of) the Laplace transform and satisfy the  
differential equations (we remind that $\go $ is primitive, $\p_i$
are flat and $\p_N$ is the unit element in the Frobenius algebra 
$(T_{\tau}\T, \bullet)$ so that $\p_N F =1$): 
\[ \p_i I = (\p_i\bullet) \p_N I,\ \ \p_N I = -\p_{\gl} I, \ \
(\gl \p_{\gl} + E) I = (\mu -k -1/2) I .\]
The equations determine the solution unambiguously from an initial condition
(the specialization to $\tau=0$ will suffice). 
Also by definition $\p_{\gl} I^{(k)}=I^{(k+1)}$.
In terms of the generating function $\f =\sum I^{(k)} (-z)^k$ the equations
read:
\begin{equation} \label{5.3}
 \p_i \f = z^{-1}(\p_i\bullet ) \f,\ \ \p_N \f+\p_{\gl}\f=0,\ \ 
(z\p_z+\gl\p_{\gl}+E) \f = (\mu-1/2) \f .\end{equation}
The specialization $\f_0=\f(\gl, 0)$ satisfies respectively
\[ \p_i \f_0=0, \ \p_{\gl }\f_0 = -z^{-1}\f_0, \ \ 
(z\p_z+\gl \p_{\gl}) \f_0 =(\mu-1/2) \f_0.\] 
Combining this with the equations for $S_{\tau}$ 
(i.e. $\p_i S = z^{-1}\p_i\bullet S$ and $(z\p_z+E)S=\mu S-S\mu$)  
we find that 
$\f = S_{\tau}\f_0$ satisfies (\ref{5.3}). Since $\f (\gl, \tau)$ and
$S_{\tau}\f_0$ coincide at $\tau=0$ by definition, the result follows.
$\square$

\section{Stationary phase asymptotics.}

Consider the vectors fields $J_{\Gb}$ on $\T $ defined by the oscillating 
integrals $\J_{\Gb}(\tau)$ via the formula 
$(J, \p_j) = \ z\p_j\J$.
As we discussed in Section $5$, when $\Gb$ runs a basis in the appropriate 
homology group $\ZZ^N$, 
the vector fields form a fundamental solution to the system 
\begin{equation} \label{6.1}
\p_i J = z^{-1} (\p_i\bullet) J,\ \ (\p_z+(E\bullet)/z^2) J = \mu J.
\end{equation}
Now we choose $\tau$ semisimple, i.e. require the function 
$F(\cdot, \tau)$ to have $N$ non-degenerate
critical points $x_{(i)}$. We denote $u_i$ 
the corresponding critical values 
(they form a local coordinate system on $\T$ called {\em canonical})
and denote $\gD_i$ the Hessians of $F(\cdot,\tau)$ at the critical
points with respect to the primitive volume form $\go$. 
Next, we construct a basis of cycles $\Gb_1,...,\Gb_N$ as follows:  
in the levels $V_{\gl,\tau}$ varying over an infinite path from 
$\gl=u_i$ toward $\gl/z\to -\infty$ avoiding other critical values,
take a parallel family of cycles {\em vanishing as $\gl $ approaches $u_i$}
and declare their union in $\CC^m$ to be $\Gb_i$. 
In fact many details do not matter here since we are going to replace 
the oscillating integrals $\J_{\Gb_i}$ by their stationary phase asymptotics
near $u_i$. In this way we get an asymptotical fundamental
solution to the same system (\ref{6.1}). The asymptotical
solution has the form $J \sim \Psi\ R_{\tau}(z)\ \exp (U/z)$ where: 
\begin{itemize} \item $U=\operatorname{diag}(u_1,...,u_N)$, 
\item $\Psi(\tau)$ is the transition matrix
from the basis $\{ \p_j \}$ in $T_{\tau}T$ 
to the basis $\sqrt{\gD_i}\p/\p u_i$ orthonormal
with respect to the residue metric, 
\footnote{ Note that the residue metric in the canonical coordinates
assumes the form $\sum \gD_j^{-1}(du_j)^2$. Respectively the matrix
$\Psi$ satisfies the orthogonality condition 
$\sum_{a,b}\Psi^a_i(\p_a,\p_b)\Psi^b_j =\delta_{i,j}$, and therefore
$[\Psi^{-1}]^j_i=\sum_a (\p_a,\p_i)\Psi^a_j=\gD_j^{-1/2}\p_i u_j$.} and
\item 
$R_{\tau}(z)=\1+R_1z+R_2z^2+...$ 
is a formal power series with matrix coefficients depending on $\tau$.  
\end{itemize}     
According to \cite{GiS} (Proposition, part (d)) an asymptotical solution
of this form to the system (\ref{6.1}) is unique and automatically satisfies
the symplectic condition $R^t(-z)R(z)=\1$. According to the definition
of the total descendent potential (\ref{1.1}) given in \cite{GiQ}
the data $\Psi, R, U$ in (\ref{1.1}) come from this unique
asymptotical solution and thus coincide with the corresponding 
ingredients of the stationary phase asymptotics 
$J\sim \Psi R(z) \exp (U/z)$ described above. 

The coefficient $C$ in the formula (\ref{1.1}) is defined 
(uniquely up to a non-zero constant factor) in terms of the
diagonal entries of the matrix $R_1$ (see \cite{GiQ}):
\[ C(\tau) := 
\exp \left( \frac{1}{2} \int^{\tau} \sum R_1^{ii}(u)du_i \right) .\]  
The genus-$1$ Gromov-Witten potential ${\mathbf F}^{(1)}$ of a semisimple
Frobenius structure mentioned in the previous section is defined 
(up to an additive constant) by
\[ {\mathbf F}^{(1)}(\tau ) := \frac{1}{48}\sum_{i} \ln \Delta_i(\tau) + 
\ln C(\tau) .\]
As it is shown, for instance, in \cite{GiE}, the function 
${\mathbf F}^{(1)}$ is constant in the case of $A_2$-singularity. 
Using Hartogs' principle one can derive from this (see, for example, \cite{H})
that for arbitrary singularity it extends analytically from semisimple points $\tau$
through the caustic. In particular, ${\mathbf F}^{(1)}$ is
constant (as a regular function on $\T$ of
zero homogeneity degree) in the case of all simple singularities.

To complete the description of the formula (\ref{1.1}), we note that
$(\q_1,...,\q_N) = \Psi^{-1} \q \in \CC^N [z]$ 
is the coordinate expression for $\q \in H [z]$ in terms of our orthonormal
basis in $T_{\tau}\T$ identified with $H=T_0\T$ via the flat metric 
$(\cdot,\cdot )$.

\medskip

We have therefore the ancestor potential defined by the formula
\begin{equation} \label{6.1'} \A_{\tau} (\q)=   
\Psi (\tau)\ \hat{R}_{\tau} e^{(U/z)\hat{\ }} 
\ \prod_{i=1}^N \D_{A_1} (\q_i) \Delta_i^{-1/48}(\tau), \end{equation}
and our next goal is to learn how to commute the vertex operators 
$\gG^{\beta}_{\tau}$ past $\Psi R e^{U/z}$.  

In fact, the conjugation $\hat{J} (e^{\f})\hat{\ }\hat{J}^{-1}$ of an 
element of the Heisenberg group by a quantized symplectic transformation
is proportional to $(e^{J^{-1}\f})\hat{\ }$. We postpone the discussion
of the proportionality coefficient and compute $J^{-1}\f$. 

Let $\gb_i$ be the cycle in $H^{2l}(V_{\gl,\tau})$ vanishing as $\gl \to u_i$
along the same path as the one participating in the definition of the
non-compact cycle $\Gb_i$. The vector $J_{\Gb_i}$ of oscillating integrals is
expressed via $I^{(l)}_{\beta_i}$ by the ``Laplace transform'' along the path:
\[ J_{\Gb_i}(\tau) = \frac{(-z)^{-l}}{\sqrt{- 2\pi  z}}  
\int_{u_i}^{-\infty} e^{\gl/z} I^{(-l)}_{\beta_i}(\gl,\tau) d\gl = 
\frac{1}{\sqrt{-2\pi  z}} 
\int_{u_i}^{\infty} e^{\gl/z} I^{(0)}_{\beta_i}(\gl,\tau) \]
(note that $I^{(k)}(u_i, \tau)=0$ for $k<0$). 
Near the critical value $\gl=u_i$ we have the expansion 
\[ (I_{\beta_i}^{(0)},\p_j) = \frac{\p_j u_i}{\sqrt{\gD_i}} 
\frac{2}{\sqrt{2(\gl-u_i)}} (1+...) \]
where the dots mean power series in $2(\gl-u_i)$. In components, we find 
\begin{equation} \label{6.2} [I^{(0)}_{\beta_i}]^j = 
\sum_{a}\Psi^j_{a} \left( \delta^{a i} + \sum_{k>0}A_k^{a i} [2(\gl-u_i)]^k 
\right) \ \frac{2}{\sqrt{2(\gl-u_i)}}.
\end{equation}
Using the change of variables $\gl-u_i=-z x^2/2$ we compute
\begin{align} \frac{2}{\sqrt{-2\pi z}} \int_0^{-\infty} e^{\gl/z} 
[2(\gl-u_i)]^{k-1/2}d\gl & =  
\frac{(-z)^{k+1/2}}{\sqrt{-2\pi z}} \ e^{u_i/z} 
\int_{-\infty}^{\infty} e^{-x^2/2} x^{2k} dx \\ 
& = (-z)^k\ (2k-1)!!\ e^{u_i/z}.\end{align}
Thus the asymptotics of $J_{\Gb_i}$ assumes the form
\[ [J_{\Gb_i}]^j \sim \sum_{a} \Psi^j_{a}
\left( \delta^{a i}+\sum_{k>0} (2k-1)!!\ A_k^{a i} (-z)^k \right) \ e^{u_i/z} ,\]
and therefore $R_k^{\a i} = (-1)^k (2k-1)!!\ A_k^{\a i}$. Substituting this
into (\ref{6.2}) and combining with the Taylor formula 
$ e^{u/z} \sum_{k\in \ZZ} z^k I^{(k)}(\gl ) = 
\sum_{k\in \ZZ} z^k I^{(k)} (\gl+u) $,
we arrive at the following result.

\medskip

{\bf Theorem 3.} {\em Near $\gl=u_i$ we have 
$\sum_{k\in \ZZ} (-z)^k I^{(k)}_{\beta_i}  = \ \Psi \ R (z) \ e^{U/z} \ 
\1_i \ {\mathbf I}$ where $\1_i =\sqrt{\gD_i} \p/\p u_i$ is the 
$i$-th unit coordinate vector in $\CC^N$ and ${\mathbf I}(z,\gl):=
2 \sum_{k\in \ZZ} (-z)^k (\frac{d}{d\gl})^k \ (2\gl)^{-1/2}$.}

\medskip

{\em Remark.} Note that coefficients of a $z$-series of the form $R \f$, where 
$R=\sum_{l\geq 0} R_l z^l$ and $\f = \sum I^{(k)} (-z)^k$, are infinite sums
$\sum_{l\geq 0} \pm R_l I^{(k-l)}$. They converge in the $\sqrt{\gl-u_i}$-adic sense
as long as $I^{(k)}$ expands near $\gl =u_i$ into a Laurent series in $\sqrt{\gl-u_i}$
such that the lowest exponent tends to $\infty$ as $k \to -\infty$.     

\section{The phase factors} 

Let us introduce {\em the phase $1$-form}
\begin{equation} \label{7.1}
\tilde{\W}_{\beta}(\gl,\tau) :=-( I^{(0)}_{\beta}(\gl,\tau), d I^{(-1)}_{\beta}(\gl,\tau) ) 
=  \sum_{i=1}^N ( I^{(0)}_{\beta}, \p_i\bullet I^{(0)}_{\beta} ) \ dt_i. 
\end{equation} 
It depends quadraticly on the cycle $\b$, and we will occasionally denote 
$\tilde{\W}_{\a,\b} = -(I^{(0)}_{\a},d I^{(-1)}_{\b})$
its polarization which is symmetric and bilinear in $\a,\b$. 
The phase form is, generally speaking, multiple-valued and is ramified along the
discriminant where $\gl$ is a critical value of $F(\cdot, \tau)$.

We discuss below some basic properties of the phase form.
\footnote{T. Milanov has found an elegant description of the phase form 
in terms of the Frobenius multiplication on the {\em co}tangent bundle. 
WE refer to \cite{GiM} for details and for 
explicit formulas in terms of the root systems in the case of
ADE-singularities.}  

\medskip

{\bf 1.} Both $\tilde{\W}_{\beta}$ and the polarizations are closed since $\p_i(\p_j\bullet) =
\p_j(\p_i\bullet)$.  

{\bf 2.} The phase form is invariant under $\p_{\gl}+\p_{N}$, i. e. $\tilde{\W }$ 
is determined 
by the restriction $\W(\tau):=\tilde{\W} (0,\tau)$ via $\tilde{\W} (\gl,\tau ) = 
\W (\tau-\gl \1)$.

{\bf 3.} Let $E=\sum \deg (t_i) t_i \p_i$ be the Euler vector field. Then 
$i_E \W_{\a,\b} = - \lan \a,\b\ran$.
Here $\lan \a,\b\ran$ is the intersection index 
normalized in such a way that the self - intersection of a vanishing cycle equals $+2$. 
Indeed, $( a , E \bullet b)$ is known to be proportional to the intersection form 
carried over to the cotangent spaces $T^*_{\tau}\T$ by the differential 
of the period map $\tau \mapsto [ d^{-1}\go ]$ defined by the primitive form $\go$
(see, for instance, \cite{H}). According to \cite{Va}, the proportionality coefficient 
is independent of the singularity and can be computed in an example.   

The property of $\W$ means that $\exp \int \W_{\beta}$ is homogeneous of degree 
$-\lan \beta,\beta \ran$. For example, when $\a$ is a $1$-point cycle in the level
$x^n/n=-\tau_{n-1}$ of the $A_{n-1}$-singularity, we have 
\begin{equation} \label{7.2}
\int_{\tau_{n-1}=-\1}^{-\gl\1} \W_{\a} = \int_{x^n=n}^{n\gl} 
\sum_{i=1}^{n-1}\frac{x^{i-1}}{x^{n-1}}\frac{x^{n-1-i}}{x^{n-1}} d(-\frac{x^n}{n}) =
\frac{1-n}{n}\ln \gl .\end{equation}
Note that $(n-1)/n$ is the self-intersection index of $\a$ projected to the {\em reduced}
homology group. 

{\bf 4.} 
Suppose that a cycle $\a$ is invariant under the monodromy along a loop in the
complement to the discriminant. Then the phase form $\W_{\a}$ is single-valued along the 
loop, and we can talk about the {\em period} $\oint \W_{\a}$. When a small loop $\c=\b^2$ 
goes twice around the discriminant near a non-singular point, then the monodromy is 
trivial, and
\begin{equation} \label{7.3} \oint \W_{\a} =-2\pi i \lan \a,\b \ran^2, \end{equation}
where $\b$ is the cycle vanishing at the corresponding critical point. 

Indeed, $\a = \lan \a,\b\ran \b/2 + \a'$, where $\lan \a',\b\ran =0$. Let $\gl=u$ be
the critical value. Then $I^{(0)}_{\a'}$ is analytic at $\gl=u_i$, and $I^{(0)}_{\b}$  
expands in $\sqrt{\gl-u}$ as in (\ref{6.2}). This implies that $\oint \W_{\a',\b}$
and $\oint \W_{\a',\a'}$ vanish, while $\oint \W_{\b/2,\b/2} = -2\pi i$ (as in (\ref{7.2}) 
with $n=2$).     

Obviously, the same is true for any conjugation $\d \b^2 \d^{-1}$ (which itself is  
the square of $\d \b \d^{-1}$).
   
\medskip

{\bf Proposition 1.} {\em In the case of a simple singularity, suppose that a cycle
$\a$ has integer intersection indices with vanishing cycles and is invariant under
the monodromy along some loop $\c$. Then the corresponding period $\oint \W_{\a}$ is an integer 
multiple of $2\pi i$.}

\medskip

{\em Proof.} If a transformation from a finite reflection group preserves some vector,
then it can be written as a composition of reflections in hyperplanes containing the 
vector. On the other hand, the (monodromy) reflection group of a simple singularity 
is known to coincide with the quotient of corresponding Artin's braid group 
(i. e. the fundamental group of the complement to the discriminant) 
by the normal subgroup generated by the squares of standard generators. Thus the loop $\c$
can be written as the composition $ \c  = \b_1^2  ... \b_r^2  \b'_1 ... \b'_s$,
where $\b_i, \b_i'$ are ``small'' loops around non-singular points of the discriminant,
and the monodromy along $b_i'$ preserves $\a$. The loops 
$\b_i' $ have zero contributions to the period $\int_{\c}\W_{\a}$
(since $\a$ is orthogonal to the corresponding vanishing cycles), while the periods
$\int_{\d^{-1}\b_i^2 \d }\W_{\a} = \int_{\b_i^2}\W_{\a} = -2\pi i \lan\a,\b_i\ran^2 $
are integer multiples of $2\pi i$. $\square $   
    
\medskip

We will show now how central constants in various commutation relations between
vertex operators and symplectic transformations are expressed in terms of the phase form.

In the situation of Theorem $1$, let us compute the factor $e^{W(\f_{+},\f_{+})/2}$
defined by the formulas (\ref{5.1}, \ref{5.2}).  
Differentiating (\ref{5.1}) and using $\p_i S (z)= z^{-1}S(z)$ and 
$(\p_i\bullet)^t=\p_i\bullet$ we find $\p_i W(\q,\q) = 
([S\q]_0,\p_i\bullet [S\q]_0)$ where $[S\q]_0$ denotes the zero mode
in $S(z)\q(z)$. Since $S|_{\tau=0}=\1$, we
see from (\ref{5.1}) that $W|_{\tau=0}=0$ and conclude
\[ W(\q,\q) = \int_0^{\tau} \sum ( [S\q]_0,\p_i \bullet [S\q]_0) dt_i.\]
The differential $1$-form here is closed and the integral does not
depend on the path connecting the origin $0\in \T$ with $\tau\in \T$
(at least when $\q \in H[z]$).
    
We apply the formula to $\q = \sum_{k\geq 0} (-z)^k I_{\beta}^{(k)}(\gl, 0)$.
According to Theorem $2$, $[S_{\tau}\q ]_0 = I_{\beta}^{(0)}(\gl, \tau)$
and therefore the exponent $W(\f_{+},\f_{+})$ in (\ref{5.2}) can be written as 
\begin{equation} \label{7.4} 
 W(\f_{+},\f_{+})=\int_0^{\tau} \sum_i (I_{\beta}^{(0)}(\gl, t), \p_i\bullet 
I_{\beta}^{(0)}(\gl, t) ) dt_i = \int_{0}^{\tau} \tilde{\W}_{\beta} \end{equation}
This integral may depend on the path (in the complement of the discriminant) 
which determines the branch of the multiple-valued vector-function $I^{(0)}_{\beta}$. 
Slightly abusing notation, we indicate the end-points in such integrals but suppress
the name of the path. However we always assume that in different integrals the path
is the same whenever the end-points are the same. Also we choose $-\1\in \T$,
(defined by $F (x, -\1) = f(x)-1$) for the base point. 

Rewriting the integral via $\W$ 
\[  W(\f_{+},\f_{+}) = \int_{-\1}^{\tau-\gl\1} \W_{\b} - \int_{-\1}^{-\gl\1} \W_{\b},\]
computing the second integral as
\[ \int_{-\1}^{-\gl\1} \W_{\b} = -\int_{1}^{\gl} (I^{(0)}_{\b}(\xi,0),I^{(0)}_{\b}(\xi,0))
d\xi = - \lan \b,\b \ran \int_{1}^{\gl} \frac{d\xi}{\xi},\]
and combining this with Theorem $2$, we arrive at the following conclusions.

{\bf Proposition 2.} {\em Introduce the vertex operator 
\begin{equation} \label{7.5} 
\gG^{\beta}_{\tau}(\gl) =e^{ \sum_{k<0} I^{(k)}_{\beta}(\gl, \tau) (-z)^k}
e^{ \sum_{k\geq 0} I^{(k)}_{\beta}(\gl, \tau) (-z)^k}.\end{equation} 
Then we have
\[ \hat{S}_{\tau}\ e^{-\lan \beta,\beta\ran \int_{1}^{\gl} d\xi/2\xi} 
\gG^{\beta}_0(\gl)\ 
\hat{S}_{\tau}^{-1} = e^{\int_{-\1}^{\tau-\gl\1}\W_{\beta}/2} \ 
\gG^{\beta}_{\tau}(\gl) .\]}

\medskip

The weights $\gl^{(1-n)/n}$ in the formulation of Theorem $1$ coincide with 
$\gl^{-\lan \a,\a\ran}$ for the $1$-point cycles $\a$ and differ from 
$\exp (-\lan\a,\a\ran \int_{1}^{\gl} d\xi/\xi )$ by the corresponding $n$-th roots of unity 
(as explained before the formulation of Theorem $1$).

\medskip

{\bf Corollary.} {\em  In the case of $A_{n-1}$-singularities 
an element $\D $ of the Fock space satisfies the
$n$KdV hierarchy (\ref{4.3}) if and only if for some --- and then for
all --- $\tau \in \T$ the corresponding elements
$\A_{\tau} =\hat{S}_{\tau}\D$ satisfy the condition
\begin{equation} \label{7.6} 
\sum_{\a} (\gG^{-\a}_{\tau} \A_{\tau})(\q') \ (\gG^{\a}_{\tau} \A_{\tau})(\q'')\ 
e^{\int_{-\1}^{\tau-\gl\1}\W_{\a} - \lan \a,\a\ran \int_{1}^{\gl} d\xi/\xi} \frac{d\gl}{\gl^{\lan \a,\a\ran}}
\ \ \ \text{is regular in $\gl$}.\end{equation} } 

\medskip

{\em Remark.} As was explained in Section $3$, the regularity condition 
refers to expansions into Laurent series
in $\gl^{-1}$, and in particular the multiple-valued functions 
$(I_{\a}^{(k)}, [\phi_i])$ and $\int \W_{\a}$ should be understood as series
expansions $\gl^{\mu_i-1/2-k} (a_0+a_1\gl^{-1}+...)$ and respectively
$-\lan \a,\a\ran \ln \gl + b_1\gl^{-1}+b_2\gl^{-2}+...$ near $\gl = \infty$.

\medskip

Let us return now to the situation of Theorem $3$.

According to \cite{GiQ}, Proposition $7.3$, 
the action of the operator $\hat{R}^{-1}$ on elements of the Fock space 
is given by the formula
\[ (\hat{R}^{-1} \G)(\q) = (e^{\h V(\p,\p)/2} \G)(R \q ) \]
where $(R\q)(z)=R(z)\q(z)$, and the ``Laplacian'' $V(\p,\p)=
\sum (\p_{q_k}, V_{kl}\p_{q_l})$ is defined by 
\[ \sum_{k,l\geq 0} V_{kl}w^kz^l = \frac{\1-R(w)R^t(z)}{w+z}.\] 
This easily implies 
\[ \hat{R}^{-1} (e^{\f})\hat{\ }\hat{R}= e^{V\f^2_{-}/2}
\ (e^{R^{-1}\f})\hat{\ },\]
where $\f_{-}=\sum_{k\geq 0} (-1)^{-1-k}(f_{-1-k}, q_k)$ 
is interpreted as a linear function of $\q$. 
When $\f  =\sum_{k\in \ZZ} I_{\beta}^{(k)}(-z)^k$, we have 
$\f_{-}=\sum (I^{(-1-k)}_{\beta},q_k)$. Using  
$\p_{\gl} I_{\beta}^{(-1-k)}=I_{\beta}^{(-k)}$, we find
\footnote{We slightly abuse notation by identifying $T_{\tau}\T$ with 
$\CC^N$ by $\Psi^{-1}$ and denoting in the same way $(\cdot,\cdot)$ 
the metric on $T_{\tau}\T$, and the standard inner products on $\CC^N$
and $\CC^{N*}$.} 
\begin{align} \p_{\gl} V \f^2_{-} & =  \sum_{k,l\geq 0} 
((I_{\beta}^{(-k)},\cdot), [V_{k-1,l}+V_{k,l-1}] (I_{\beta}^{(-l)},\cdot)) \\ 
 & = (I_{\beta}^{(0)},I_{\beta}^{(0)}) -
(\sum R^t_k I_{\beta}^{(-k)}, \sum R^t_l I_{\beta}^{(-l)} ). \end{align}
Let us assume now that $\beta$ is a vanishing cycle $\beta_i$. By Theorem $3$,
\[ \sum R^t_k I_{\beta_i}^{(-k)} = \sum R^t_k \sum (-1)^l 
R_l (\frac{d}{d\gl})^{-l-k} \frac{2\ \1_i}{\sqrt{2(\gl-u_i)}} =
\frac{2\ \1_i}{\sqrt{2(\gl-u_i)}} .\]
Also $V\f_{-}^2 =0 $ at $\gl=u_i$ since 
$\f_{-k} \sim (\gl-u_i)^{k+1/2}(\1_i + ...)$ vanish at $\gl=u_i$.
Thus
\[ V \f_{-}^2 = \int_{u_i}^{\gl} \left(  
(I^{(0)}_{\beta_i}(\xi,\tau), I^{(0)}_{\beta_i}(\xi,\tau))
- \frac{2}{(\xi - u_i)}\right)\ d\xi .\]   
Note that near $\xi=u_i$ both integrals diverge, but in the same way,
so that the difference converges. 
The integral can be rewritten as 
\[ V\f_{-}^2 = \int_{\tau-\gl\1}^{\tau-u_i\1} \left( 
\W_{\beta_i} - \frac{2 dt_N}{\tau_N-u_i-t_N} \right) .\]
Finally, conjugation of vertex operators by $e^{(U/z)\hat{\ }}$ 
is a special case of Proposition $2$ and has the following effect: 
$ \1_i /\sqrt{2(\gl-u_i)}$ is transformed to $\1_i /\sqrt{2\gl}$, and the 
corresponding factor $e^{W(\f_{+},\f_{+})/2}$ is equal to $(\gl-u_i)/\gl$.
Note that the correspondence between the branches of the $\sqrt{\cdot }$ depends
on the choice of a path connecting $\gl-u_i$ with $\gl$. 

We summarize.

\medskip

{\bf Proposition 3.} {\em Let $\b_i$ be one of the vanishing cycles. Put 
\begin{equation} \label{7.7}
 W_i := \int_{\tau-\gl\1}^{\tau-u_i\1} \left( \W_{\beta_i/2} - 
\frac{dt_N}{2(\tau_N-u_i-t_N)} \right), \ \ \ w_i=\int_{\gl-u_i}^{\gl}\frac{d\xi}{2\xi}
\end{equation}
Then
\[ \left(\Psi \hat{R} e^{(U/z)\hat{\ }}\right)^{-1}\ 
e^{-W_i/2}\ \gG^{\pm_{\beta_i/2}}_{\tau} \ 
\left(\Psi \hat{R} e^{(U/z)\hat{\ }}\right) =  e^{-w_i/2} \  
\left( ... \1 \otimes (\gG^{\pm})_{(i)} \otimes 
\1 ... \right),\] 
where $\gG^{\pm}$ are the vertex operators (\ref{3.7}), and the subscript $(i)$
indicates the $i$-th position in the tensor product. }
  
\medskip

{\em Remark.} The integration path in the definition of $w_i$ is the same as the
one that determines the branch of $\sqrt{\cdot }$ under the translation $\sqrt{\gl-u_i} 
\mapsto \sqrt{\gl}$.

\medskip

Now let us consider a cycle $\a$ represented as the sum $c\ \b/2 + \a'$ where 
$\b$ is the cycle vanishing over the point $(\gl,\tau)=(u,\tau)$ on the discriminant,
and $\a'$ is a cycle invariant under the local monodromy near this point (so that 
$c=\lan \a,\b\ran$).

\medskip

{\bf Proposition 4.} {\em For the vertex operators (\ref{7.5}) we have
\[ \Gamma_{\tau}^{\a} = e^{c \int_{\tau-\gl\1}^{\tau-u\1} \W_{\b/2,\a'}}
 \Gamma_{\tau}^{\a'} \Gamma_{\tau}^{c \b/2} .\]}

\medskip

{\em Proof.} It is clear that $\Gamma_{\tau}^{\a} =e^{K} 
\Gamma_{\tau}^{\a'}\Gamma_{\tau}^{c \b/2}$. The proportionality coefficient $e^{K}$
arises from commuting $e^{\hat{\f}_{-}}$ across $e^{\hat{\g}_{+}}$, where
$\f =c\ \sum (-z)^k I^{(k)}_{\b/2}$ and $\g =\sum (-z)^l I^{(l)}_{\a'}$. The constant
$K$ is equal therefore to the symplectic inner product $\gO (\f_{-},\g_{+})$. One 
easily finds  
\[ K = c\ \sum_{k\geq 0} (-1)^k (I_{\b/2}^{(-1-k)}, I^{(k)}_{\a'}).\]
On the other hand, consecutive integration by parts yields
\[ \int_u^{\gl} (I^{(0)}_{\beta/2}, I^{(0)}_{\a'} )\ d\xi =
\sum_{k=0}^{m-1}  (-1)^k(I^{(-1-k)}_{\beta/2},I^{(k)}_{\a'})|_u^{\gl} 
+(-1)^m  \int_u^{\gl} (I^{(-m)}_{\beta/2}, I^{(m)}_{\a'} )\ d\xi.\]
Note that $I^{(-1-k)}_{\beta/2} \sim (\gl-u)^{k+1/2}(\1_i+...)$
and vanish at $\gl=u$, while $I^{(k)}_{\a'}$ are holomorphic at $\gl=u$.
Thus the last integral is $o (\gl-u)^{m-1/2}$ and hence tends to $0$ 
as $m\to \infty$. We conclude that 
\[ K=c\ \int_u^{\gl} (I^{(0)}_{\b/2}(\xi,\tau), I^{(0)}_{\a'}(\xi,\tau)) \ d\xi
= - c\ \int_{\tau-u\1}^{\tau-\gl\1} \W_{\b/2,\a'} .\] 

\section{Asymptotical elements of the Fock space}

Various expressions with quantized symplectic transformations and vertex operators contain
numerous infinite sums, and we have to discuss now precise meaning of our formulas.

By an {\em asymptotical function} we will mean an expression of the form
\[ \exp \sum_{g\geq 0} \h^{g-1} \F^{(g)}(\t ),  \]
where $\F^{(g)}$ is a formal function on the space $H[ \t ]$ of polynomials 
$\t (z) =t_0+t_1z+t_2z^2+...$ with vector coefficients 
$t_k=\sum_{\a} t_k^{\a}\phi_{\a} \in H$. 

We will say that an asymptotical function is {\em tame} if
\[ \frac{\p}{\p t_{k_1}^{\a_1}} ... \frac{\p}{\p t_{k_r}^{\a_r}} |_{\t =\0} \F^{(g)} = 0
\ \ \text{whenever} \ \ k_1+...+k_r > 3g-3+r .\]
In particular, each $\F^{(g)}$ is a formal series $\sum F^{(g)}_{a,b} (t_0)^{a}(t_1)^{b}$
of $t_0,t_1$ with the coefficients which are {\em polynomials} on $t_2,...,t_{3g-2+|a|}$.

The Witten -- Kontsevich tau-function is tame (as well
as ancestor potentials \cite{GiQ} in Gromov -- Witten theory are --- because 
$\dim_{\CC} \M_{g,r} =3g-3+r$). 

An asymptotical function is identified with {\em an asymptotical element in the Fock space}
(in the formalism of Section $2$) via the {\em dilaton shift} $\q(z)=\t(z)-z$ and becomes
therefore an asymptotical function of $\q $ (tame or not)
with respect to the {\em shifted origin} $\q =-z$. The notation $-z:=(-\1)z$ is the only 
place where we use that the space $H$ contains a distinguished non-zero vector $\1 $.

\medskip

{\bf Proposition 5.} {\em Let $R$ be an {\em upper-triangular} element of the 
twisted loop group, i. e. $R(z)=\1+R_1z+R_2z^2+...$, and $R^t(-z)R(z)=\1$.
Then the action of the quantized operator $\hat{R}$ on tame asymptotical elements of 
the Fock space is well-defined and yields tame asymptotical elements.}

\medskip
   
{\em Proof.} As mentioned in Section $7$ the action of $\hat{R}^{-1}$ on an asymptotical 
function $\G$ takes the form
\[ (\hat{R}^{-1} \G) (\t) = (e^{\h V(\p,\p)/2} \G) (R\t + \c),\ \ \text{where} \ \ 
\c(z)=z-R(z)z.\]
The operation $\ln G \mapsto \ln ( e^{\h V(\p,\p)/2} \G)$ can be described in terms
of summation over connected graphs with vertex contributions defined by partial derivatives
of $\ln G:= \sum \h^{g-1}\F^{(g)}$, and edge factors given by the coefficients $V_{kl}$ of 
the ``Laplacian'' $V(\p,\p)$.
\footnote{We are not going to enter here a detailed discussion of the Wick formula
underlying the graph summation technique. However the reader may track the origin
of the ``graphical'' interpretation of the operator $R$ back to \cite{GiS}.}  
In order to check that $\ln (e^{\h V(\p,\p)/2} \G)$ is tame,
let us examine the contribution of a connected graph with $E$ edges 
into a Taylor coefficient at $t^{\a_1}_{k_1}...t^{\a_r}_{k_r}$. 
Let 
\begin{itemize}
\item $g(v)$ be the genus of a vertex $v$, $e(v)$ --- the number of edges incident 
to the vertex ($\sum e(v) =2 E)$, 
\item $l(v)$ --- the total sum of the indices in the derivatives 
$V_{kl}\p_{t_k}\p_{t_l}$ applied to the vertex $v$ ($\sum l(v)=: L$), 
\item $r(v)$ --- the number of marked points in $v$ ($\sum r(v) =r$), 
\item $k(v)$ --- the total sum of the indices among $k_1,...,k_r$ attributed to the vertex
($\sum k(v)=k_1+...+k_r=: K$).
\end{itemize} 
The total genus $g$ of the graph (i.e. the power 
of $\h$ to which the graph contributes) is determined by the formula 
$g-1=\sum (g(v)-1) + E$. We see that $g\geq 0$ since $g(v)\geq 0$ and 
$E - \sum_v 1 \geq -pr0file5
1$. Since $\G$ is tame, the contribution of 
the graph vanishes unless $k(v)\leq 3g(v)-3+e(v)+r(v)-l(v)$ for each $v$. Summing up we 
find 
\[ K\leq 3\sum (g(v)-1)+2E+r = 3g-3 + r - L - E \leq 3g-3+r.\]
Thus the required condition is satisfied. Moreover, the number of
edges of the graph and the indices in the edge factors $V_{k,l}$ are bounded 
($L+E\leq 3g-3+r$). Thus $\ln (e^{\h V(\p,\p)/2} \G)$ is well-defined since there
are only finitely many terms of each genus $g$ and degree $r$.   

The substitution of $R(z) \t(z)$ instead of $\t(z)$ preserves the above conclusions
since the multiplication by $R=R_0+R_1z+R_2z^2+...$ does not decrease the indices 
$k_1,...,k_r$ (determined by the degree in $z$).

Finally, the series $z-R(z)z$ starts with $z^2$ since $R_0=\1$. Therefore the dilaton 
shift $\t(z) \mapsto \t(z) + z-R(z)z$ is also a well-defined operation in the class of
tame asymptotical functions. $\square$  
 
\medskip

As it is mentioned in Section $5$, {\em lower-triangular} operators $\hat{S}^{-1}_{\tau}$ 
act on an asymptotical element $\G$ of the Fock space by $e^{W(\q,\q)/\h} \G ([S\q]_{+})$.
The change  $\bar{\q} = [S\q]_{+}$ means $\bar{\t} (z) = [S_{\tau}(z) \t(z)]_{+} -\tau$
or, in components, $\bar{t}_0 =\sum S_k(\tau) t_k-\tau$, $\bar{t}_1 =\sum S_k t_{k+1}$, ...
Suppose that $\ln \G$ is a formal function of $\bar{\t}$
and is therefore defined in the formal neighborhood of $\bar{\t}=0$.
When $\t(z)$ is a polynomial, $\bar{\t} = S_{\tau}\t -\tau =\0$ means 
$t_0=\tau, t_1=t_2=... =0$. This makes $\ln (\hat{S}^{-1}_{\tau}\G )$ 
a well-defined formal function of $t_0-\tau, t_1,t_2, ...$ (and $\h$).

The operators $\hat{S}$ with $S(z)=S_0+S_1z^{-1}+S_2z^{-2}+...$ do {\em not} preserve 
the class of tame functions. In particular this applies to the rightmost operator
in $\Psi \hat{R} \exp (U/z)\hat{\ }$. Yet the formula (\ref{6.1'}) for the ancestor 
potential makes sense and defines a tame asymptotical function $\A_{\tau}$ because 
the operators $\exp (u/z)\hat{\ }$ {\em preserve} $\D_{A_1}$.
Indeed, the string equation for the Witten -- Kontsevich tau-function coincides with
$(1/z)\hat{\ } \D_{A_1} = 0$.

More generally, let us call a tame asymptotical function $\G$ 
{\em $T$-stable}, if $T\G$ is also tame. Let $\G$ be 
$\exp (U/z)\hat{\ } $-stable for all diagonal matrices $U$.
\footnote{These requirements are satisfied, for example, if
$\G=\D_1 (\q_1)...\D_N(\q_N)$ where $\D_i$ are obtained from $\D_{A_1}$ by translations
$\q \mapsto \q + \a$, where $\a(z)=a_0+a_1z+a_2z^2+...$ is a vector-polynomial 
(or even a series) with coefficients which are formal $\h$-series such that $a_0$ and $a_1$
are smaller than $1$ in the $\h$-adic norm (and $a_k \to 0$ in this norm as 
$k \to \infty$).} 
Then 
$e^{{\mathbf F}^{(1)}}\A_{\tau}:=\Psi_{\tau} \hat{R_{\tau}} \exp (U(\tau)/z)\hat{\ } \G $ 
are well-defined and tame, while 
$\hat{S}_{\tau}^{-1} \A_{\tau}$ are defined as asymptotical functions of $\t (z) - \tau$.
Moreover, according to Theorem $7.1$ in \cite{GiQ}, the asymptotical element
$\D:= \prod \Delta_i^{-1/48} \hat{S}_{\tau}^{-1}\A_{\tau}$ does not depend on $\tau$ 
and is therefore well-defined as an asymptotical function of $(t_0,t_1,...)$ in
the formal neighborhood of $(\tau,0,...)$ with any {\em semisimple} $\tau$.

\medskip

Let us examine now the regularity condition (see Corollary to Proposition $2$ of 
Section $7$) in the description of integrable hierarchies via vertex operators.
The action of the vertex operators of the form 
$\Gamma_{\tau}^{-\beta}\otimes \Gamma_{\tau}^{\beta}$ 
on functions $\G (\x') \otimes \G (\x'')$
is described more explicitly (see (\ref{3.1}), (\ref{3.2})) 
as composition of translations and multiplications:
\begin{equation} \label{8.1}
\exp\left(2 \sum_{k\geq 0} (I^{(-1-k)}_{\beta}, \frac{q_k}{\sqrt{\h}})\right) 
 \exp\left( -\sum_{k\geq 0} (-1)^k I^{(k)}_{\beta} \sqrt{\h} \p_{q_k} \right)  
\G (\x + \q) \G (\x -\q) . \end{equation}    
The coefficient $I^{(k)}_{\beta}(\gl,\tau)$ can be represented near $\gl =\infty$
by an infinite series in fractional powers $\gl^{\nu}$ with the exponents $\nu$ from the
union of $N$ arithmetical sequences $\mu_i-1/2-k+\ZZ_{-}$. As we remarked in Section $8$,
the phase factors $\exp \int \W_{\beta}$ also expand into such series with 
$\nu \in -\lan \beta,\beta \ran +\ZZ_{-}$. The formulation that a vertex operator
expression like (\ref{7.6}) is {\em regular in $\gl$}
instructs us to expand (\ref{8.1}) into a $\q$-series. In fact (\ref{8.1}) is 
manifestly invariant under the classical monodromy operator. As a result, 
the coefficient at a given monomial $\q^{\m}$ expands into a Laurent series in $\gl^{-1}$ 
(since the coefficient at each power $\gl^{\nu}$ depends only on finitely many 
$I^{(k)}_{\beta}$). The regularity condition, by definition, means that the
coefficients at negative powers of $\gl$ vanish (so that the Laurent series in 
$\gl^{-1}$ is a polynomial in $\gl$). 
     
On the other hand, recalling the genus expansion 
$\G =\exp (\sum \h^{g-1}\F^{(g)})$ and using the notation $Q_k:=q_k/\sqrt{\h}$,
we can rewrite (\ref{8.1}) as
\begin{equation} \label{8.2}
\exp \left[ 2\sum_{k\geq 0} (I^{(-1-k)}_{\beta}, Q_k)   
+ \sum_{g\geq 0} \h^{g-1} \sum_{\pm} \F^{(g)} \left( \x \pm \sqrt{\h} \Q 
\mp\sqrt{\h} \sum_{k\geq 0} I^{(k)}_{\beta} (-z)^k  \right) \right]
\end{equation}
The functions $\F^{(g)}$ are formal series of $\x$. Rewriting the exponent as a series in $\h$
we see that the $\h^{-1}$-term $2\F^{(0)}(\x)$ does not depend on $\gl$ 
and all the $\h^{-1/2}$-terms cancel out. 

\medskip

{\bf Proposition $6$.} {\em Suppose that $\G=\exp \sum \h^{g-1} \F^{(g)}$ is a tame 
asymptotical function of $\x$. Then (\ref{8.1}) divided by $\exp (2\F^{(0)}(\x)/\h)$
expands into a power series in $\sqrt{\h}, \x$ and $\Q$ whose coefficients depend
polynomially on finitely many $I^{(k)}_{\beta}$ each.}

\medskip  

{\em Proof.} Recall that expansions of each $\F^{(g)}(\x)$ with as power series in $x_0,x_1$
have coefficients which depend only on finitely many $x_2,x_3,...$. Note that each 
$I^{(0)}_{\beta}, I^{(1)}_{\beta}$ in (\ref{8.2}) brings with itself an extra $\sqrt{\h}$. 
We conclude that modulo high powers of $\sqrt{\h}$ the exponent of (\ref{8.2}) 
is a series in $\Q, \x$ whose coefficients depend polynomially on finitely many 
$I^{(k)}_{\beta}$ each. Subtracting the singular term $2\F^{(0)}(\x)/\h $ 
and exponentiating does not alter this conclusion. $\square $

\medskip

Proposition $6$ means, that the regularity requirement, when applied to tame asymptotical 
functions, can be understood not only as 
a statement about expansions near $\gl =\infty$, but also as the property
of analytic functions of $\gl$ (the polynomial expressions of $I^{(k)}_{\beta}$
and of the phase factors $\exp \int \W_{\beta}$) to be single-valued polynomial
functions of $\gl$. 
    
\medskip

Finally, it is worth reiterating here some of our remarks from Sections $5$ and $6$
about conjugations of vertex operators by quantized elements of the twisted loop
group:

\begin{itemize}

\item the conjugation 
$\hat{S}_{\tau}^{-1} \Gamma_{\tau}^{\beta} \hat{S}_{\tau}$ by lower-triangular
elements is well-defined via the expansion of
$\sum_k I^{(k)}_{\beta}(\gl,\tau) (-z)^k$ as a series near $\gl=\infty$, 

\item the conjugation $\hat{R}^{-1}_{\tau} \Gamma^{\beta_i}_{\tau} \hat{R}_{\tau}$
by upper-triangular elements is well-defined in terms of expansions 
near the critical value $\gl=u_i$, and

\item the conjugation by $\exp (u_i/z)\hat{\ }$ acts on $\sum_k I^{(k)}(\gl) (-z)^k$ 
as the translation $\gl \mapsto \gl+u_i$; it is applied in our computations only 
to the vertex operator defined by the analytic functions 
$I^{(k)}(\gl) = (d/d\gl)^k (\gl-u_i)^{-1/2}$.

\end{itemize}

\section{From $n$KdV to $n-1$ KdV}

We prove here Theorem $1$ as a special case (with $\D_1=...=\D_{n-1}=\D_{A_1}$)
of a more general result which yields a solution of the $n$KdV hierarchy from $n-1$ 
solutions of the KdV hierarchy.

\medskip

{\bf Theorem $4$.} {\em Suppose that asymptotical functions $\D_i(\q_i), i=1,...,n-1$,
are tame and stable with respect to the string flows
$e^{(u_i/z)\hat{\ }}$. Let us assume that the ingredients $C, S, \Psi, R$ and $U$ of the
formula (\ref{1.1}) correspond to the Frobenius structure of the $A_{n-1}$-singularity. Then 
\begin{equation} \label{9.1} 
\D := C(\tau) \hat{S}_{\tau}^{-1} \Psi (\tau) \hat{R}_{\tau} e^{(U(\tau)/z)\hat{\ }}
\prod_{i=1}^{n-1} \D_i(\q_i) \end{equation}
satisfies the equations of the $n$KdV-hierarchy: 
\begin{equation} \label{9.2}
\left[ \sum_{\text{$1$-point cycle $\a$}} \Gamma^{-\a}_{0}\otimes \Gamma^{\a}_{0} \ 
\gl^{-\lan \a, \a\ran} d\gl \right] \ \D \otimes \D \ \ \text{is regular in $\gl$}. 
\end{equation} } 

\medskip

{\em Proof.} Similarly to Corollary from Theorem $2$ and Proposition $2$, it is sufficient to
prove that $\A_{\tau} := e^{-{\mathbf F}^{(1)}} \hat{S}_{\tau} \D$ satisfies the condition:
\begin{equation} \label{9.3} 
\left[ \sum_{\text{$1$-point cycle $\a$}} \Gamma^{-\a}_{\tau}\otimes \Gamma^{\a}_{\tau} \ 
e^{\int_{-\1}^{\tau-\gl\1}\W_{\a}+\lan\a,\a\ran\int_1^{\gl} \frac{d\xi}{\xi}} \gl^{-\lan \a,\a\ran}d\gl \right] \ 
\A_{\tau} \otimes \A_{\tau} \ \ \text{is regular in $\gl$} \end{equation}  
for at least one value of $\tau$. We choose $\tau $ to be generic (so that $F(\cdot, \tau)$
is a Morse function) and prove (\ref{9.3}) as follows. 

In view of Proposition $6$ we can interpret (\ref{9.3}) in terms of analytic functions in
$\gl$ (rather than series in $1/\gl$). Since all the $n$ one-point cycles $\a$ form an orbit 
of the monodromy group of the $A_{n-1}$-singularity, one can argue that (\ref{9.3}) is 
invariant under the whole monodromy group (and not only the classical monodromy operator). 
Thus (\ref{9.3}) is meromorphic with possible poles at the distinct critical values 
$u_1,...,u_{n-1}$. The regularity property will follow if we prove that there are no poles 
at $\gl=u_i$.

Let $\beta = \a_{+}-\a_{-}$ be the cycle vanishing at $\gl=u_i$, and
$\a_{\pm}$ are two of the $n$ one-point cycles. If $\a \neq \a_{\pm}$,
then $\a$ is invariant under the monodromy around $u_i$, the
corresponding vector-functions $I^{(k)}_{\a}$ are holomorphic at $\gl=u_i$,
and therefore the phase factor and respectively the whole summand in (\ref{9.3}) 
with the index $\a$ is holomorphic at $\gl=u_i$ as well.  

When $\a=\a_{\pm}$, we have $\a=\pm \beta/2 + \a'$ where 
$\a'=(\a_{+}+\a_{-})/2$ is invariant under the monodromy around $u_i$. 
Thus $I^{(k)}_{\a_{\pm}} = I^{(k)}_{\pm \beta/2} + I^{(k)}_{\a'}$
where the second summand is holomorphic at $\gl=u_i$. We have therefore
\begin{equation} \label{9.4} \gG^{\a_{\pm}}_{\tau} =e^{\pm K}\ 
\gG^{\a'}_{\tau}\ \gG^{\pm \beta/2}_{\tau} ,\end{equation}
where the proportionality coefficient $e^{\pm K}$ is described by Proposition $4$ 
(with $u=u_i$, and $c=\pm 1$). 
Thus the two summands in (\ref{9.3}) with $\a =\a_{\pm}$ add up to
\begin{equation} \label{9.5} 
\gG^{-\a'}_{\tau}\otimes \gG^{\a'}_{\tau}\left[ \left( \sum_{\pm} C_{\pm}(\gl)
\gG^{\mp \beta/2}_{\tau}\otimes \gG^{\pm \beta/2}_{\tau} \right) 
\A_{\tau} \otimes \A_{\tau} \right]\ d\gl , \end{equation}
where $C_{\pm}$ are some phase factors combined from (\ref{9.3}) and (\ref{9.4}).

Let us now recall that 
\[ \A_{\tau} = \Psi \hat{R} e^{(U/z)\hat{\ }} \prod \D_i(\q_i) \D_i^{-1/48} \]
and apply Theorem $3$. We see that 
\begin{itemize}
\item the square bracket in (\ref{9.5}) has the form of the operator
$\Psi \hat{R} e^{(U/z)\hat{\ }}$ applied to a product 
$\prod \F_i$ of $n-1$ functions in $n-1$ different 
groups of variables $(\q'_i,\q''_i)$, 
\item the factors $\F_i$ corresponding to $i$ with $\beta_i\neq \beta$ 
are equal to $\D_i(\q'_i)\D_i(\q''_i)$,
\item the factor corresponding to $\beta_i=\beta$ has the form  
\begin{equation} \label{9.6}  \sum_{\sqrt{2\gl}=\pm \zeta} c_{\pm} (\gl) 
(\gG^{-} \D_i)(\q'_i)\ (\gG^{+}\D_i)(\q''_i)\ \frac{d\gl}{\sqrt{\gl}},  \end{equation}
\item the phase factors $c_{\pm}/\sqrt{\gl}$ come from $C_{\pm}$ and from the phase
factors described by Proposition $3$.
\end{itemize}

We assume that the factors $\sqrt{\gl}$ here differ by the sign (rather than coincide)
the same way as in (\ref{3.6}) (or (\ref{4.3}) when $n=2$).  

\medskip

{\em We claim that near $\gl=u_i$ the functions $c_{\pm}(\gl)$ coincide, are
single-valued and analytic.}

In order to justify the claim, let us compute the phase factors explicitly.
We have: 
\begin{align} \ln c_{\pm} = 
 \int_{-\1}^{\tau-\gl\1} \W_{\a_{\pm}} + \lan\a_{\pm},\a_{\pm}\ran \int_1^{\gl}
\frac{d\xi}{\xi} - \ln \gl^{\lan \a_{\pm},\a_{\pm} \ran} + \ln \sqrt{\gl}-
\int_1^{\gl}\frac{d\xi}{2\xi} \\ 
 \pm 2\int_{\tau-\gl \1}^{\tau-u_i\1} \W_{\b_i/2,\a'}
+ \int_{\tau-\gl\1}^{\tau-u_i\1} ( \W_{\b_i/2} -\frac{dt_N}{2(\tau_N-u_i-t_N)} ) 
-\int_{\gl-u_i}^{1} \frac{d\xi}{2\xi}.\end{align}

Using bi-linearity of the phase form $\W$ with respect to the cycles 
$\a_{\pm}=\pm \b_i/2 +\a'$ we rewrite:
\begin{align} \label{9.7} \ln c_{\pm}  = 
\int_{-\1}^{\tau-\gl \1} \W_{\a'} &   + \lan \a_{\pm},\a_{\pm}\ran 
\oint_{\c_{\pm}}\frac{d\xi}{\xi} + \oint_{\c'_{\pm}} \frac{d\xi}{2\xi}  \\ 
\label{9.8}  \pm 2 \int_{-\1}^{\tau-u_i\1}  \W_{\b_i/2,\a'} 
+  \int_{-\1}^{\tau-(u_i+1)\1}   & \W_{\b_i/2} + 
\int_{\tau-(u_i+1)\1}^{\tau-u_i\1} (\W_{\b_i/2}- \frac{dt_N}{2(\tau_N-u_i-t_N)} ) \\
 \label{9.9}   & +  \int_{\tau_N-(u_i+1)}^{\tau_N-\gl}\frac{dt_N}{2(\tau_N-u_i-t_N)} 
- \int_{\gl-u_i}^1\frac{d\xi}{2\xi} . \end{align} 
The constant $u_i+1$ is chosen to make the integrals in (\ref{9.9}) cancel exactly.
The contours $\c_{\pm},\c'_{\pm}$ in (\ref{9.7}) as well as all terms in 
(\ref{9.8}) may depend on the cycle $\a_{\pm}$ but are independent of 
$\gl$, while the first integral in (\ref{9.7}) is a function of $\gl $ analytic at $\gl=u_i$ 
and independent of the cycle. 
This implies that the phase factors $c_{\pm}$ are proportional to each other and are
analytic near $\gl=u_i$. 

\medskip

Let us show that the proportionality coefficient equals $1$. Since $\b_i/2=(\a_{+}-\a_{-})/2$
and $\a'=(\a_{+}+\a_{-})/2$, we have $4 \W_{\b_i/2,\a'}=\W_{\a_{+}}-\W_{\a_{-}}$ and 
therefore
\begin{equation} \label{9.10} 
\ln c_{+}-\ln c_{-} =  \lan\a_{-},\a_{-} \ran \oint_{\c_{+}-\c_{-}} \frac{d\xi}{\xi} + 
\int_{-\1}^{\tau-u_i\1}( \W_{\a_{+}}-\W_{\a_{-}}) + \oint_{\c'_{+}-\c'_{-}} \frac{d\xi}{2\xi}.
\end{equation} 
Note that the one-point cycles $\a_{\pm}$ belong to the same orbit of the classical
monodromy (i.e. the cyclic group of the Coxeter transformation) and therefore the first 
integral in (\ref{9.10}) can be interpreted as $\oint_{\c_1} \W_{\a_{-}}$ where the loop 
$\c_1$ makes several turns about $\gl=0$ inside the line $-\gl \1$ so that $\a_{-}$ 
transported along the loop becomes $\a_{+}$ in the end. 

Let $\c_2(\ge)$ denote the path starting at $-\1$ and approaching the point $\tau-u_i\1$ on the 
discriminant (as in the second term in (\ref{9.10})) but stopping a small distance $\ge $ away 
from it. Let $\c_3(\ge)$ be a loop of size $\ge$ going around the discriminant near $\tau-u_i\1$
(so that $\a_{+}$ transported along $\c_3$ becomes $\a_{-}$ in the end).
The integral $\int \W_{\a_{+}}$ along the path $\c_2(\ge)\c_3(\ge)\c_2^{-1}(\ge)$ does not
depend on $\ge$ (for homotopy reasons). In the limit $\ge \to 0$ it spits out the middle term
of (\ref{9.10}) plus $\lim_{\ge\to 0}\int_{\c_3(\ge)} \W_{\a_{+}}$. Writing 
$ \W_{\a_{+}} = \W_{\b_i/2} + 2\W_{\b_i/2,\a'}+\W_{\a'}$ near $\tau-u_i\1$ we see that
the first summand contains the term $d\gl /2(\gl-u_i)$, and the rest is 
either analytic at $\gl=u_i$ or has a  
singularity like $(\text{analytic function}) \times d\gl /\sqrt{\gl-u_i}$. This implies
that 
\[ \lim_{\ge \to 0} \int_{\c_3(\ge)} \W_{\a_{+}} = \oint \frac{d\gl}{2(\gl-u_i)} = \pi \sqrt{-1},\]
which coincides with the last integral in (\ref{9.10}). 
We conclude that (\ref{9.10}) can be interpreted as the period of $\W_{\a_{-}}$ 
along the loop $\c_1 \c_2 \c_3 \c_2^{-1}$.
The cycle $\a_{-}$ is invariant under the monodromy along this loop.
According to Proposition $1$ the period is an integer multiple of $2\pi \sqrt{-1}$.
Thus $c_{+}=c_{-}$. 

\medskip

The proof of  Theorem $4$ is now completed as follows. Since $\D_i$ satisfy the KdV 
hierarchy, we conclude that (\ref{9.6}) is regular in $\gl$. 
This implies that $\prod \F_i$, and
hence (\ref{9.5}) is single-valued near $\gl=u_i$ and has no pole at $\gl=u_i$. 
Since the other ingredients of (\ref{9.3}) are also holomorphic at $\gl=u_i$, 
we find that (\ref{9.3}) is regular at $\gl=u_i$. In particular, (\ref{9.3}) is invariant
with respect to the whole monodromy group (regardless of reliability of the previously
mentioned abstract argument) and is regular in $\gl$. $\square$  

\section{Some applications}

Due to Theorem $1$ the total descendent potential $\D_{A_{n-1}}$ defined by (\ref{1.1}) 
satisfies the $n$KdV-hierarchy and is therefore ``a tau-function''. In addition it 
satisfies the string equation $(1/z)\hat{\ } \D_{A_{n-1}}=0$ (due to \cite{GiQ}).
Solutions of the $n$KdV-hierarchies satisfying the string equation have been studied 
in the literature (see for instance \cite{W, Sch}) under the name {\em $W_n$-gravity}. 
By definition, the tau-functions in the $W_n$-gravity theory are formal functions
of the variables $t_0,t_1,t_2, ... \in H$. Our functions $\D_{A_{n-1}}$, to the contrary,
are known to expand in formal series near {\em semisimple} $t_0$. 
It is our present goal to identify $\D_{A_n}$ with the tau-function singled out in the
theory of $W_n$-gravity, and in particular  --- to establish analyticity of 
the total descendent potential at $t_0=0$. 

\medskip

It will be convenient for us to use another form of the $n$KdV-hierarchy based on the 
concept of Baker functions. Let $\exp \sum_{g\geq 0} \h^{g-1}\F^{(g)}(\t)$ 
be an asymptotical function in a formal neighborhood of $\t ={\mathbf 0}$. 
Given an asymptotical function $\G $, the corresponding {\em Baker function } 
\cite{SW} (or {\em wave function} \cite{K}) is defined as 
\[ b_{\G} = (\Gamma^{\a} \G)/\G = e^{-\sum_{k\geq 0} (I^{(-1-k)}_{\a}, q_k)/\sqrt{\h}} 
\G (\q + \sqrt{\h} \sum_{k\geq 0} I_{\a}^{(k)}(-z)^k ) /\G(\q)   .\]
Here $\Gamma^{\a}$ is the vertex operator (\ref{4.1}) corresponding to a one-point
cycle $\a$. The Baker function can be understood as a $q$-series
\[ b_{\G} =\sum b_{\G}^{(\mathbf m)}\q^{\mathbf m} \]
with coefficients $b_{\G}^{({\mathbf m})}$ which are Laurent series of 
$\zeta^{-1}=\gl^{-1/n}$ 
(whose coefficients, in their turn, are Laurent series in $\sqrt{\h}$).
Let $\CC_{\sqrt{\h}}((\zeta^{-1}))$ be the space of all such Laurent series, and
let $V_{\G}$ denote the subspace spanned over $\CC_{\sqrt{\h}}[\gl]$ by the coefficients 
$b_{\G}^{({\mathbf m})}$. 
According to the grassmannian description \cite{SW} of the KP-hierarchy, 
$\G$ satisfies the $n$KdV-hierarchy if and only if $V_{\G}$ belongs to the principal 
cell of the semi-infinite grassmannian (i.e. projects isomorphically onto 
$\CC_{\sqrt{\h}}[\zeta]$ along $\zeta^{-1}\CC_{\sqrt{\h}}[[\zeta^{-1}]]$).
  
On the other hand, conjugation of $\Gamma^{\a}(\gl) \gl^{-\lan \a,\a \ran/2}$ 
by the string flow $\exp (u/z)\hat{\ }$ yields $\Gamma^{\a}(\gl-u) (\gl-u)^
{-\lan\a,\a\ran/2}$ and therefore 
\begin{equation} \label{10.1} [ (1/z)\hat{\ }, \Gamma^{\a} ] = 
d/d\gl -\lan\a,\a\ran/2 \gl .\end{equation}
In particular the string flow $\exp (u/z)\hat{\ }$ acts on the vertex operator 
expression in (\ref{4.3}) by translation $\gl \mapsto \gl-u$ and therefore preserves the 
regularity requirement in (\ref{4.3}). Thus the string flow is a symmetry of the
$n$KdV-hierarchy. Moreover, let $\G_{\tau}$ be the total descendent potential 
$\D_{A_{n-1}}$ considered as an asymptotical function in the formal variable 
$\t =  \q -\tau + z$, where $\tau =\sum_{i=1}^{n-1} \tau_i [\phi_i] \in H$ 
is a semisimple point (and $z$ represents the dilaton shift). 
The invariance of $\D_{A_{n-1}}$ with respect to the string flow  
can be restated via the Baker function $b_{\G_{\tau}}$ as invariance
of the space $V_{\G_{\tau}}$ with respect to the operator  
\[
 A:=\frac{d}{d\gl} - \frac{\lan \a,\a\ran }{2\gl} + \frac{(\tau, I^{(0)}_{\a}(\gl))}
{\sqrt{\h}} - \frac{(\1, I^{(-1)}_{\a}(\gl))}{\sqrt{\h}} .\]
The first two terms here come from (\ref{10.1}) and the others come from 
\[  e^{\p_{\tau}-z\p_{\1}} \Gamma^{\a} = \exp \{ \h^{-1/2} [(\tau, I^{(-1)}_{\a})-
(\1, I_{\a}^{(-2)})]\}\ \Gamma^{\a}e^{\p_{\tau}-z\p_{\1}} .\]
Note that the space $V_{\G_{\tau}}$ contains $v:=b_{\G_{\tau}}|_{\q={\mathbf 0}}$ 
which is a power series in $\zeta^{-1}$ with the constant term $1$. Such series form a group
acting on the semi-infinite grassmannian via multiplication. The invariance of $V_{\G_{\tau}}$
with respect to $A$ is equivalent to invariance of 
$U:=v^{-1}V_{\G_{\tau}}$ relative to $B=v^{-1} A v$. We have
\begin{equation} \label{10.2} 
B=\frac{d}{d\gl} + \frac{(\1, I^{(-1)}_{\a}(\gl))}{\sqrt{\h}}
 + \frac{(\tau, I^{(0)}_{\a}(\gl))}{\sqrt{\h}} -\frac{\lan \a,\a\ran }{2\gl}+ 
\sum_{k\geq 0} (-1)^k \frac{(f_k,I_{\a}^{(k+1)})}{\sqrt{\h}}, \end{equation}
where the linear function $\sum_{k\geq 0} (-1)^k (f_k,q_k)$ of $\q$ is the differential of
$\h \ln \G_{\tau}$ at $\t =0$.  
The space $U$ is a free $\CC_{\sqrt{\h}}[\gl]$-module of rank $n$. It contains the series $1$ 
and hence contains all $\gl^m$ and all $B^k(1)$. 
Note that the $\zeta^{-1}$-series $B^k(1)$ starts with $\zeta^k$ 
(since $(\1,I^{(-1)}_{\a})\sim \gl^{1/n}$) and therefore 
$1, B(1),...,B^{n-1}(1)$ form a basis in $U$, while 
$B^n(1) =\h^{-n/2} n \gl + \sum_{k>0} a_k \gl^{1-k/n} $. 
The coefficients $a_1,...,a_{n}$ are uniquely determined by $\tau$. 
Representing $B^n(1)-\h^{-n/2} n \gl $ as a linear combination of the basis vectors,
we obtain a system of equations for $f_0,f_1,...$. It is straightforward to see that the
system is triangular and unambiguously determines all $f_k$ via $\tau$. We find 
that the space $U$ and respectively $V_{\G_{\tau}}$ is unique for each $\tau$.
Due to the correspondence 
between semi-infinite subspaces, Baker functions and tau-functions 
(see for instance \cite{SW} or Exercises $14.44$ -- $14.47$ in \cite{K})
we conclude that {\em the asymptotic function $\G_{\tau}$ is completely 
characterized up to a scalar factor as a formal solution to the $n$KdV-hierarchy 
near $\q=\tau-z$ satisfying the string equation.} 

Let us consider now the tau-function function $\G_{\tau}$ corresponding to the 
(non-semisimple) $\tau=0$. Existence of the function and of the corresponding space 
$V_{\G_0}$ follows from the results of \cite{Sch} (or from the above argument which
is a slight variation on the theme of \cite{Sch} anyway). Note that the corresponding
operator
\begin{equation} \label{10.3} A = d/d\gl - \h^{-1/2}(n\gl)^{1/n}-(n-1)/2n\gl \end{equation}
is homogeneous (of degree $-1$) with respect to the grading 
$\deg \gl =1,\ \deg \h = 2+2/n$.
This implies that the basis $A^k(v), \ k=0,...,n-1$, in $V_{\G_0}$ and 
respectively the tau-function $\G_0$ are homogeneous in the appropriate sense.
More explicitly, $\ln \G_0$ has the form $\sum \h^{g-1} \F^{(g)}(\t )$, 
where $\F^{(g)}$ are formal series of $t_{i,k}, \ i=1,...,n-1,\ k=0,1,2, ...$ 
homogeneous of degree $(1-g)(2+2/n)$ with respect to the grading 
$\deg t_{i,k} = (i+1)/n\ -\ k$. This follows from the famous fact \cite{SW} that
the flows of the KP-hierarchy (which in our notation are represented by the derivations
$\sqrt{\h} \p_{t_{i,k}}$) correspond in the grassmannian description to the multiplication
by $\zeta^{kn+n-i}$ (and hence $(kn+n-i)\deg \zeta = \deg \sqrt{\h} - \deg t_{i,k}$).

By definition, the asymptotical function $\G_0$ is ``the tau-function of 
the $W_n$-gravity theory'' and, according to a conjecture of E. Witten \cite{W}, coincides
with the total descendent potential in the intersection theory (developed in \cite{PV})
on moduli spaces of complex curves equipped with $n$-spin structures. 

Consider now the formal homogeneous function $\h \ln \G_0$ as a power series in 
$\h, t_1,t_2,...$ with coefficients (which are therefore also homogeneous)
depending on $t_0=(t_{i,0},...,t_{n-1,0})$ . Since all the components
of $t_0$ have positive degrees, we conclude that each coefficient is {\em polynomial}
in $t_0$. Thus translations $\G_0 (t_0+\tau, t_1,t_2,...)$ are well-defined and yield
asymptotical functions satisfying the same conditions --- 
the $n$KdV-hierarchy and the string equation --- as $\G_{\tau}(t_0,t_1,t_2,...)$.
The previous uniqueness argument now implies $\G_{\tau}(\t) = \G_0(\t+\tau)$ for 
all $\tau\in H$. We have proved the following result.

\medskip

{\bf Theorem 5.} {\em The total descendent potential $\D_{A_{n-1}}$ of the 
$A_{n-1}$-singularity coincides with the tau-function $\G_0$ introduced in the
$W_n$-gravity theory.}

\medskip

{\bf Corollaries.} {\em (1) The total descendent potential $\D_{A_{n-1}}$ of the
$A_{n-1}$-singularity (which is an asymptotical function of $t_0,t_1,...$
defined in a formal neighborhood of $(t_0,0,...)$ with semisimple $t_0$)  
extends across the caustic to arbitrary $t_0 \in H$.

(2) The ancestor potentials $\A_{\tau}=\hat{S}_{\tau}\D_{A_{n-1}}$ are
well-defined for all $\tau\in H$.

(3) The descendent potential $\D_{\A_{n-1}} =\A_{0}$ and is tame.

(4) The Gromov -- Witten potentials ${\mathbf F}^{(g)}$ of the $A_{n-1}$-singularity
(defined by $\sum \h^{g-1} {\mathbf F}^{(g)} (\tau) := \ln \A_{\tau} |_{\t={\mathbf 0}}$)
are polynomial functions of $\tau \in H$ of weighted degree $(1-g)(2+2/n)$ and 
therefore vanish for $g>1$.} 
  
%
%

\section*{Appendix: Dispersionless limit}

In the ``dispersionless limit'' $\h \to 0$
Theorem $1$ implies that the genus-$0$ descendent potential
$\F^{(0)}$ of Saito's Frobenius structure on the miniversal deformation of the 
$A_{n-1}$-singularity satisfies the dispersionless $n$KdV-hierarchy.  
We give here a more direct proof of this fact using only the general 
theory of $n$KdV-hierarchies and the results of Section $5$. 
No doubt, this relationship between the dispersionless nKdV hierarchies 
and $A_{n-1}$-singularities has been known for quite a while (see 
for instance, \cite{D, Kr}), but we are not so sure about the  
following lemmas. 

\medskip
 
Let us recall from Section $3$ that an asymptotical function 
$\Phi (\x) = e^{\sum \h^{g-1} \phi^{(g)}(\x)}$ is said to satisfy the KP-hierarchy if
\begin{equation} \label{a1} 
\operatorname{Res}_{\zeta =\infty} d\zeta \ 
e^{2\sum_{j>0}\zeta^jy_j/\sqrt{\h}} \ 
e^{-\sum_{j>0}\frac{\zeta^{-j}}{j} \sqrt{\h} \p_{y_j}}
\ \Phi(\x+\y)\Phi(\x-\y) =0. \end{equation} 

\medskip

{\bf Lemma 1.} {\em A function $\phi^{(0)}$ satisfies the dispersionless limit of the 
KP-hierarchy ($n$KdV-hierarchy) if and only if for each $\q$ the function 
$\exp [ (d^2_{\q}\phi^{(0)})(\x)/2\h ] $,
where $d^2_{\q}\phi$ is the quadratic form of the $2$nd differential of $\phi^{(0)}$ at $\q$,
satisfies the KP-hierarchy (the $n$KdV-hierarchy respectively).}

\medskip

{\em Proof.} In order to pass to the limit $\h\to 0$, divide (\ref{a1}) by $\Phi^2(\x)$, put 
$Y:=\y/\sqrt{\h}$ (and respectively $\sqrt{\h}\p_{\y}=\p_{Y}$) and expand
\[ \frac{\Phi (\x+\sqrt{\h}Y)}{\Phi(\x)}\frac{\Phi(\x-\sqrt{\h}Y)}{\Phi(\x)} = 
e^{ W(Y) + O (\h) }, \]
where $W(Y)$ is the quadratic form $d^2_{\x}\phi^{(0)}$. Taking $\h=0$ results in a closed
system of equations for $\phi^{(0)}$ which, by definition, is the {\em dispersionless} 
KP-hierarchy. Namely, $\phi(\x)$ satisfies the differential equations of the 
dispersionless hierarchy if for all $\x$ the quadratic differential 
$W=d^2_{\x}\phi$ satisfies the system of algebraic equations  
\begin{equation} \label{a2} 
\operatorname{Res}_{\zeta =\infty} d\zeta \ 
e^{2\sum_{j>0}\zeta^jY_j} \ 
e^{-\sum_{j>0}\frac{\zeta^{-j}}{j} \p_{Y_j}}
\ e^{W (Y)} =0. \end{equation}      
It is an observation of T. Milanov that the condition (\ref{a2}) for a quadratic form $W$
is equivalent to (\ref{a1}) for the corresponding {\em Gaussian distribution}  
$\Phi (\x) = e^{W(\x)/2\h}$. (One may in fact take $\h=1$ everywhere.) 

Solutions of the (dispersionless) $n$KdV-hierarchy are those solutions of the
(dispersionless) KP-hierarchy which do not depend on $x_i$ with $i$ divisible by $n$. 
$\square $

\medskip

{\em Remark.} The form (\ref{a1}) of the KP-hierarchy is stronger than the usual system 
of dynamical equations for the function $u:=(\ln \Phi)_{xx}$ (here $x=x_1$). For example, the 
KdV-hierarchy in the form $u_{x_i}=(\L_i(u,u_x,u_xx,...))_x$ is automatically satisfied by
any $\Phi = \exp W/2$ since $u$ is constant. We will see in a moment
that this is not at all the case for the algebraic system (\ref{a2}).

\medskip  

Which Gaussian distributions $e^{\sum_{ij}W_{ij}x_ix_j/2}$ (we put $\h=1$) 
satisfy the KP- and $n$KdV-hierarchies?
Consider the corresponding Baker function 
\[ b_{W}(\x) := e^{\sum \zeta^j x_j} e^{-\sum \zeta^{-j}\p_{x_j}/j} e^{W(\x)/2}
= b_W({\mathbf 0}) e^{ \sum_i x_i(\zeta^i-\sum_j W_{ij}\zeta^{-j}/j) },\]
where $b_W({\mathbf 0}) = \exp (\sum_{ij} W_{ij}\zeta^{-i-j}/2ij)$.
Let $V_W$ denote the subspace in $\CC ((\zeta^{-1}))$ spanned by the Taylor coefficients
of the {\em normalized} Baker function $b_W(\x)/b_W({\mathbf 0})$.

\medskip

{\bf Lemma 2.} {\em A Gaussian distribution $\exp W/2 $ satisfies the KP-hierarchy
($n$KdV-hierarchy)
if and only if the corresponding normalized Baker function generates a semi-infinite
subspace $V_{W}$ which is a subring (respectively a $\CC [\zeta^n]$-subalgebra) 
in $\CC ((\zeta^{-1}))$.}

\medskip

{\em Proof.} Indeed, the subspace $V_W$ being semi-infinite (which is necessary and 
sufficient for a function to satisfy the KP-hierarchy) means that the Laurent series
\[ 1, \ \ \zeta-\sum W_{1j}\zeta^{-j}/j, \ \ \zeta^2-\sum W_{2j}\zeta^{-j}/j,\ \ ... \]
form a basis in $V_{W}$. Taylor coefficients of the normalized Baker
function are arbitrary products of these series which therefore have to be in $V_{W}$.

Solutions to $n$KdV-hierarchy correspond to semi-infinite subspaces invariant with
respect to multiplication by $\zeta^n$.  $\square$.
      
\medskip

Note that 

(i) the above basis in the space $V_W$ is canonical in the sense that it is obtained by
lifting the basis $1,\zeta,\zeta^2,...$ from $\CC [\zeta]$ to $V_{W}$ along 
$\zeta^{-1}\CC [[ \zeta^{-1}]]$,   

(ii) the Baker function of a space $V$ from the principal cell of the semi-infinite 
grassmannian is normalized {\em iff} the $1$st element in the canonical basis is $1$,

(iii) the rest of the basis determines the coefficients $W_{ij}$ unambiguously,
which establishes a $1$--$1$ correspondence between Gaussian distributions satisfying 
the KP-hierarchy and semi-infinite subrings $V \subset \CC ((\zeta^{-1}))$ from the principal
cell of the grassmannian,

(iv) $V_{W} =\CC [x]$, where $x:= \zeta - \sum W_{1j}\zeta^{-j}/j$, and $W_{1j}, j=1,2,...$,
are arbitrary numbers.

\medskip

{\bf Corollary.} {\em Gaussian distributions satisfying the $n$KdV-hierarchy are
in $1$--$1$ correspondence with equations of the form 
\begin{equation} \label{a3} 
x^n+\tau_1x^{n-2}+...+\tau_{n-1} = \gl, \end{equation}
parameterized by $\tau=(\tau_1,...,\tau_{n-1})$.}

\medskip

{\em Proof.} When $V_{W}=\CC [x]$ corresponds to a solution of the $n$KdV-hierarchy,
we must have $\zeta^n\in \CC[x]$ and therefore 
$\zeta^n=x^n+\tau_0x^{n-1}+...+\tau_{n-1}$ for some $\tau_0,...,\tau_{n-1}$. On the other
hand, $W_{1j}$ must vanish for all $j$ divisible by $n$. Since all $n$ 
solutions to the
equation have the form $x(\epsilon \zeta)$ where $\epsilon$ runs through the
$n$th roots of $1$, 
we conclude that the sum $-\tau_0$ of all the $n$ solutions vanishes.

Vice versa, solving the equation for $x$ by perturbation theory near 
$x|_{\tau=0}=\gl^{1/n}$ yields a series $x=\zeta+\sum_{j\geq 0} w_j(\tau) \zeta^{-j}$
in $\zeta=\gl^{1/n}$. Since the sum of all the $n$ solutions $x(\epsilon \zeta)$ 
equals $0$, we conclude that $w_j=0$ for all $j$ divisible by $n$. 
The semi-infinite subspace $\CC[x(\zeta)] \subset \CC ((\zeta))$ 
is invariant under the multiplication by $\zeta^n=\gl$ due to (\ref{a3}). $\square$   

\medskip

The genus-$0$ descendent potential $\F^{(0)}$ of a Frobenius manifold 
(constructed in \cite{D}) can be described (due to Proposition $5.3$ and 
Corollary $5.4$ in \cite{GiQ}) in terms of the function $W$ discussed in Section $7$:
\begin{equation} \label{a4} 
W_{\tau}(\q,\q) = \int_0^{\tau} \sum ([ S_t\q]_0, \p_i\bullet [S_t\q]_0) dt_i 
\end{equation}
Namely, let us regard $W/2$ as a family of functions in $\tau \in H$ depending (quadraticly)
on the parameter $\q \in \H_{+}=H[z]$. Then $\F^{(0)}$ is the {\em critical value function} 
for this family. More precisely, the critical points $\tau$ are given by the equations
$([S_{\tau}\q]_0, \p_i\bullet [S_{\tau}\q]_0) = 0$ for all $i$. This is equivalent to
$[S_{\tau}\q]_0\bullet [S_{\tau}\q]_0 =0$ and is satisfied whenever $[S_{\tau}\q]_0=0$. 
Recall that $[S\q]_0=q_0+S_1q_1+S_2q_2+...$ where
$S=\1+S_1z^{-1}+...$, $\q =q_0+q_1z+...$. When $\q (z) = t_0-z$, we have $[S_{\tau}\q]_0
=t_0-\tau$ and find a critical point $\tau=t_0$. In general the equation $[S_{\tau}\q]_0=0$
has a unique solution $\tau (\t )$ defined by perturbation theory as a formal function of $\t=\q+z$ 
(dilaton shift). Then $\F^{(0)}(\t)=W_{\tau (\t)}/2$.

In fact the quadratic differential $d^2_{\t}\F^{(0)}$ coincides with the quadratic form
$W_{\tau(\t)}$. In particular, it depends only on the critical point $\tau$ (rather than
the parameter value $\t$).
\footnote{Moreover, according to \cite{CG, GiF}, 
Frobenius structures equipped with the genus-$0$
descendent potentials have the following axiomatic characterization. Let $\L$ denote 
the (germ at $-z$ of a) Lagrangian section in $T^*\H_{+}$ defined as the graph of $d\F^{(0)}$ 
(subject to the dilaton shift). Identifying $T^*\H_{+}$ with $(\H, \Omega)$ by means of the 
standard polarization $\H= \H_{+}\oplus \H_{-}$, we may regard $\L$ as a Lagrangian 
submanifold in $\H=H((z^{-1}))$. Then {\em $\L $ is a cone with the vertex at the origin 
and such that $\L$ intersects its tangent spaces $L$ along $zL$.} In particular, $\L$
is swept by the spaces $zL$ varying in $\dim L/zL =\dim H$-parametric family, and the
tangent spaces to $\L$ along each $zL$ are constant and coincide with $L$.} We are
going to show that in the case of $A_{n-1}$-singularities the Gaussian distributions
$e^{W_{\tau}/2}$ defined by (\ref{a4}) satisfy the $n$KdV-hierarchy --- of course, modulo the
rescaling (\ref{4.2}): $q_{i,k}=i(i+n)...(i+kn) x_{i+kn}$.

\medskip

{\bf Lemma 3.} {\em The normalized Baker function of the Gaussian distribution
$e^{W_{\tau}/2}$ corresponding to (\ref{a4}) is equal to 
\begin{equation} \label{a5} \exp \left[ \Omega \left( S_{\tau}(z) \q(z), 
\ \sum_{k\geq 0} I^{(-1-k)}_{\a}(\gl, \tau)(-z)^{-1-k}\right) \right] .\end{equation} }

\medskip

{\em Proof.} By definition, the Baker function $b_{W_{\tau}}$ is 
$e^{-W_{\tau}/2\h}\Gamma^{\a} e^{W_{\tau}/2\h}$ which after normalization and at $\h=1$ becomes
\begin{equation} \label{a6} e^{ - \sum_{k\geq 0} (I_{\a}^{(-1-k)}(\gl,0), q_k) } 
e^{ W_{\tau} (\sum_{k\geq 0} (-z)^k I^{(k)}_{\a}(\gl,0), \q )} .\end{equation}
Theorem $2$ from Section $5$ says that 
\[ S_{\tau}(z) \sum (-z)^k I^{(k)}_{\a}(\gl,0) = \sum (-z)^k I^{(k)}_{\a}(\gl,\tau) .\]
On the other hand, $d I_{\a}^{(m-1)}=-a\w I^{(m)}_{\a}$ and $d S = a \w S / z$
where $a=\sum (\p_i\bullet) dt_i =a^t$. In particular $d [S\q]_m = a\w [S\q]_{m+1}$.
Therefore computing the second exponent in (\ref{a6}) from (\ref{a4}) and
integrating by parts we find
\[ \int_0^{\tau} (I_{\a}^{(0)}(\gl,t), a(t) \w [S_t\q]_0) = 
\int_0^{\tau} (a\w I^{(0)}_{\a},[S\q]_0)
=-\int_0^{\tau}(d I^{-1}_{\a}, [S\q]_0) = \] 
\[ -(I^{-1}_{\a},[S\q]_0)|_0^{\tau}+
\int_0^{\tau}(I^{-1}_{\a},a\w [S\q]_1) = ... =
- \sum_{k\geq 0} (I^{(-1-k)}_{\a}, [S\q]_k)|_{0}^{\tau} .\]
(The integral term eventually disappears because $\q$ is a polynomial in $z$.)
The value at the lower limit $t=0$ cancels with the first exponent in (\ref{a6}), and 
the value at $t=\tau$ coincides with (\ref{a5}). $\square$

\medskip

{\bf Corollary.} {\em The vector space $V_{W_{\tau}}$ corresponding to the normalized
Baker function (\ref{a5}) is spanned by $1$ and by the components 
$(I^{(-1-k)}_{\a}(\gl,\tau),[\phi_i])$ of the period maps $I^{(m)}_{\a}$ 
with $m< 0$.} 

\medskip 

The components of $I^{(-1)}_{\a}$ are periods of the differential $0$-forms 
$x$, $x^2/2$, ..., $x^{n-1}/(n-1)$ on the level sets 
\begin{equation} \label{a7} 
\frac{x^n}{n}+\tau_1x^{n-2}+...+\tau_{n-1} =\gl \end{equation}
in the miniversal deformation of the $A_{n-1}$-singularity. 
In the case when $\a$ is a one-point cycle (i.e. $x$), 
the $\CC [\gl]$-module generated by $1,x,x^2,...,x^{n-1}$ is a subring in $\CC ((\gl^{1/n}))$
due to (\ref{a7}). It remains to show therefore that this subring coincides with $V_{W_{\tau}}$,
i.e. that it contains all components of $I_{\a}^{(m)}$ for $m<-1$. Thus the following lemma
completes the proof.

\medskip

{\bf Lemma 4.} {\em The period maps $I^{(m)}_{\a}$ satisfy the equation
\[  (\mu+1/2-m) I^{(m-1)}_{\a} = ( \gl-E \bullet )  I^{(m)}_{\a} \]
where $E=\sum_i (\deg \tau_i)\tau_i\p_{\tau_i}$ is the Euler field and $\mu+1/2$ is the
 spectral matrix, i.e. the diagonal matrix with entries $1/n,2/n,...,n-1/n$.}

\medskip

{\em Proof.} In view of the equations $\p_i I = (\p_i\bullet) \p_{n-1} I$ and 
$\p_{n-1}I=-\p_{\gl}I$ satisfied by all $I^{(k)}_{\a}$, the lemma is a reformulation of 
the homogeneity condition $(\gl \p_{\gl} + E) I^{(k)}_{\a} = (\mu-1/2-k) I^{(k)}_{\a}$
discussed in Section $5$. $\square$         

\medskip

{\em Remark.} According to a uniqueness result of Dubrovin and Zhang \cite{DZ},
the total descendent potential of a {\em semisimple} Frobenius manifold is completely 
characterized as an asymptotical function $\exp \h^{g-1} \F^{(g)}$ which satisfies 
(i) the Virasoro constraints,
(ii) the so-called $3g-2$-jet condition, and
(iii) whose genus-$0$ part $\F^{(0)}$ coincides with the genus-$0$ descendent
potential of the Frobenius manifold (constructed in \cite{D}).
According to \cite{GiQ}, the function $\D_{A_{n-1}}$ satisfies (i),(ii),(iii)
and thus would coincide with the tau-function $\G_0$ of the $W_n$-gravity theory (see Section $10$),
if $\G_{0}$ were shown to satisfy (i),(ii),(iii) as well.
In fact, the Virasoro constraints for $\G_0$ are well-known (see for instance \cite{Sch}) 
and follow from the invariance of the corresponding semi-infinite subspace 
$V_0 \subset \CC_{\sqrt{h}}((\gl^{-1/n}))$ under the operators $\gl^m A$ 
(where $A$ is given by (\ref{10.3})).
It is plausible (although at the moment we don't know a 
direct proof of this) that the $3g-2$-jet property, which is equivalent to 
the ancestor potentials $\A_{\tau}:=\hat{S}_{\tau}\G_0$ being tame for all $\tau\in H$,
can be derived from a Lax-type description of the $n$KdV-hierarchy.
Thus, since the results of this Appendix imply (iii), this would give another proof of Theorems 
$1$ and $5$. Also, Dubrovin and Zhang have informed the author that (yet another?) proof of 
these results can be obtained on the basis of their axiomatic theory of integrable 
hierarchies \cite{DZ}.

\enddocument